\documentclass[times]{TTP_DSL2006}

\renewcommand{\large}{\fontsize{14}{18pt}\selectfont}
\renewcommand{\small}{\fontsize{11}{13.6pt}\selectfont}

\newcommand{\titleformat}{\sffamily\bfseries \large}						
\newcommand{\authorformat}{\sffamily \large}							
\newcommand{\keywordsformat}{\noindent \small \sffamily}				
\newcommand{\abstractformat}{\noindent \textbf}						
\newcommand{\contentformat}{\rmfamily \normalsize\vspace{18pt}}			
\newcommand{\email}{\sffamily \small \vspace{-8pt}}						
\renewcommand{\subsection}{\textbf}	

\usepackage{url}
\usepackage{graphicx,graphics}
\usepackage{subfigure}
\usepackage{epsfig}
\usepackage{amsmath,amssymb,amsfonts}
\usepackage{mathrsfs}
\usepackage{listings}    
\usepackage{color}
\usepackage{float}
\usepackage{tabularx}
\usepackage{multirow}
\usepackage{stmaryrd}

\newcommand{\Eref}[1]{Equation (\ref{#1})}
\newcommand{\fref}[1]{Figure (\ref{#1})}
\newcommand{\Erefs}[1]{Equations (\ref{#1})}

\newcommand{\rmd}{\mathrm{d}}
\newcommand{\bveps}{\boldsymbol{\varepsilon}}
\newcommand{\bvsig}{\boldsymbol{\sigma}}
\newcommand{\pphi}{\boldsymbol{\phi}}

\newcommand{\aaa}{\mathbf{a}}
\newcommand{\bb}{\mathbf{b}}
\newcommand{\bigb}{\mathbf{B}}
\newcommand{\bc}{\mathbf{C}}
\newcommand{\dd}{\mathbf{D}}
\newcommand{\ee}{\mathbf{E}}
\newcommand{\ff}{\mathbf{F}}
\newcommand{\kk}{\mathbf{K}}
\newcommand{\bn}{\mathbf{N}}
\newcommand{\cn}{\mathbf{n}}

\newcommand{\qq}{\mathbf{q}}
\newcommand{\bq}{\mathbf{Q}}

\newcommand{\uu}{\mathbf{u}}
\newcommand{\vv}{\mathbf{v}}
\newcommand{\xx}{\mathbf{x}}

\begin{document}

\title{\titleformat Representation of singular fields without asymptotic enrichment in the extended finite element method}

\author{\authorformat Sundararajan Natarajan \inst{1}$^{,\rm{a}}$\text{,} Chongmin Song \inst{,2}$^{,\rm{b}}$\text{,} }

\institute{\sffamily School of Civil and Environmental Engineering, The University of New South Wales, Sydney, NSW 2052, Australia.}

\maketitle

\begin{center}
\email{ $^{\rm a}$sundararajan.natarajan@gmail.com, $^{\rm b}$c/song@unsw.edu.au}
\end{center}

\contentformat

\abstractformat{Abstract.} In this paper, we replace the asymptotic enrichments around the crack tip in the extended finite element method (XFEM) with the semi-analytical solution obtained by the scaled boundary finite element method (SBFEM). The proposed method does not require special numerical integration technique to compute the stiffness matrix and it improves the capability of the XFEM to model cracks in homogeneous and/or heterogeneous materials without a priori knowledge of the asymptotic solutions. A heaviside enrichment is used to represent the jump across the discontinuity surface. We call the method as the extended scaled boundary finite element method (xSBFEM). Numerical results presented for a few benchmark problems in the context of linear elastic fracture mechanics show that the proposed method yields accurate results with improved condition number. A simple MATLAB code is annexed to compute the terms in the stiffness matrix, which can easily be integrated in any existing FEM/XFEM code.

\keywordsformat{{\textbf{Keywords:}} scaled boundary finite element method, extended finite element method, boundary integration, singular fields, partition of unity methods

\vspace{-6pt}

\section{Introduction}
\vspace{-2pt}

Understanding the phenomenon of \textit{material failure} is the key to design new materials to meet the increasing demand for high strength-to-weight and high stiffness-to-weight materials. In the classical approach, the governing equations are represented by partial differential equations (PDEs). Closed form or analytical solutions to these equations are not obtainable for most problems. Numerical methods such as the finite element method (FEM)~\cite{zienkiewicztaylor2000}, the boundary element method (BEM), meshless methods, spectral methods~\cite{boyd2000}, discrete element methods~\cite{munjiza2004}, scaled boundary finite element method (SBFEM)~\cite{deekswolf2002,wolfsong2001} are essential to analyze the equations. Amongst these, the FEM is a very popular and widely used approach, because of its versatility. The FEM relies on a conforming topological map and special elements have been developed for problems involving strong discontinuities. A constant remeshing is inevitable when the discontinuities evolve in time. It has been noted that the modification (augmentation, enrichment) of the finite element spaces will improve the behaviour. An intensive research over the past 3-4 decades has led to some of the robust methods available to model crack or discontinuities, for example extended (XFEM)/ Generalized FEM (GFEM)/ Partition of Unity FEM (PUFEM)~\cite{babuvskamelenk1997,belytschkoblack1999}, $hp-$ cloud~\cite{garciafancello2000} to name a few. The additional functions, also called the '\textit{enrichment functions}' are augmented through the partition of unity framework. These additional functions carry with them the information of the local nature of the problem. Since its inception, the XFEM has been applied successfully to moving boundary problems~\cite{duddubordas2008,aubertinrethore2010,chessasmolinski2002}, to name a few. Although the XFEM is robust, it still suffers from certain difficulties and continuous research is carried out to improve the capability of the method. See Section \ref{xfemoverview} that briefly discusses the various difficulties and attempts by various researchers to overcome the difficulties. For more detailed discussion, interested readers are referred to the literature~\cite{belytschkogracie2009} and references therein.

Independently, Wolf and Song~\cite{wolfsong2001} developed a fundamental solution-less method, called the SBFEM for elasto-statics and elasto-dynamic problems. The SBFEM has emerged as an attractive alternate to model problems with singularities~\cite{songtinloi2010}.  It combines the best features of the boundary element method and the finite element method. Numerical solutions are sought on the boundary, whilst the solution along the radial lines emanating from the `\textit{scaling center}' is represented analytically. Moreover, by utilizing the special features of the scaling center, the method allows the computation of stress intensity factor directly. 

\subsection{Approach}
In this paper, we propose the xSBFEM, which combines the best features of the XFEM and the SBFEM, with the aim of improving the capability of the XFEM in representing the singular fields.
As opposed to enriching a small region in the vicinity of the crack tip with asymptotic functions (known a priori either analytically or computed numerically), in the proposed method, that region is replaced with the semi-analytical solution obtained by the scaled boundary finite element method (SBFEM). The efficiency and the accuracy of the proposed method is illustrated by solving benchmark problems taken from linear elastic fracture mechanics. The effect of the size of the SBFEM domain on the global convergence in the energy is also studied numerically. Some of the salient features of the proposed method are:

\begin{itemize}
\item Does not require a priori knowledge of the asymptotic expansions of the displacements and/or the stress fields. 
\item As it does not require asymptotic enrichments, the need for special numerical integration technique is suppressed.
\item As only boundary information is required, the total number of unknowns is reduced.
\item The scaled condition number of the stiffness matrix of the xSBFEM is comparable and/or better than the XFEM.
\end{itemize}


\subsection{Outline}
The paper is organized as follows. In the next section, we briefly recall the basic equations of the XFEM. Section~\ref{sbfem} discusses the scaled boundary finite element method (SBFEM). The proposed technique, the eXtended SBFEM (xSBFEM) is presented in Section \ref{xsbfemcou}. Two different methods for computing the stress intensity factors, other than the interaction integral are discussed in Section~\ref{sifcalculation}. The efficiency and convergence properties of the proposed method are illustrated in Section ~\ref{numexample} with a few benchmark problems taken from linear elastic fracture mechanics, followed by some concluding remarks in the last section.

\vspace{-6pt}

\section{Overview of the extended finite element method} \label{xfemoverview}
\vspace{-2pt}
Consider $\Omega \subset \mathbb{R}^2$, the reference configuration of a cracked linearly isotropic elastic body (see \fref{fig:bodywcrack}). The boundary of $\Omega$ is denoted by $\Gamma$, is partitioned into three parts $\Gamma_u, \Gamma_n$ and $\Gamma_c$, where Dirichlet condition is prescribed on $\Gamma_u$ and Neumann condition is prescribed on $\Gamma_n$ and $\Gamma_c$.  The governing equilibrium equations for a 2D elasticity problem with internal boundary, $\Gamma_c$ defined in the domain $\Omega$ and bounded by $\Gamma$ is
\begin{equation}
\nabla_s^{\rm T} \bvsig + \bb = \mathbf{0} ~\hspace{0.5cm} \in ~\hspace{0.5cm} \Omega
\label{eqn:EE}
\end{equation}
where $\nabla_s(\cdot)$ is the symmetric part of the gradient operator, $\mathbf{0}$ is a null vector, $\bvsig$ is the stress tensor and $\bb$ is the body force. The boundary conditions for this problem are:
\begin{eqnarray}
\bvsig \cdot \cn = \overline{\mathbf{t}} \hspace{0.5cm} \textup{on} ~\hspace{0.5cm} \Gamma_t \nonumber \\
\uu = \overline{\uu} \hspace{0.5cm} \textup{on} ~\hspace{0.5cm} \Gamma_u \nonumber \\
\bvsig \cdot \cn = \overline{\mathbf{t}} \hspace{0.5cm} \textup{on} ~\hspace{0.5cm}  \Gamma_c
\end{eqnarray}
where $\overline{\uu} = (\overline{u}_x,\overline{u}_y)^{\rm T}$ is the prescribed displacement vector on the essential boundary $\Gamma_u$; $\overline{\mathbf{t}} = (\overline{t}_x,\overline{t}_y)^{\rm T}$ is the prescribed traction vector on the natural boundary $\Gamma_t$ and $\cn$ is the outward normal vector. In this study, it is assumed that the displacements remain small and the strain-displacement relation is given by $\bveps = \nabla_s \uu$. Let us assume a linear elastic behaviour, the constitutive relation is given by $\bvsig = \dd \colon \bveps$, where $\dd$ is a fourth order elasticity tensor. \Eref{eqn:EE} is called the `strong form' and computational solutions of \Eref{eqn:EE} rely on a process called `discretization' that converts the problem into system of algebraic equations. The first step in transforming \Eref{eqn:EE} into a discrete problem is to reformulate \Eref{eqn:EE} into a suitable variational equation. This is done by multiplying \Eref{eqn:EE} with a test function. Let us define,
\begin{eqnarray}
\mathcal{U} = \left\{ \uu \in H^1(\Omega) \hspace{0.5cm} \textup{such that} \hspace{0.5cm}  \uu|_{\Gamma_u} = \overline{\uu} \right\} \nonumber \\
\mathcal{V} = \left\{ \vv \in H^1(\Omega) \hspace{0.5cm} \textup{such that} \hspace{0.5cm}  \vv|_{\Gamma_u} = \mathbf{0} \right\}
\end{eqnarray}
The space $\mathcal{U}$ in which we seek the solution is referred to as the `trial' space and the space $\mathcal{V}$ is called the `test' space. Let us define a bilinear form $a(\cdot , \cdot)$ and a linear form $\ell(\cdot)$,
\begin{equation}
a(\uu,\vv) \colon= \int\limits_{\Omega} \bvsig(\uu) \colon \bveps(\vv)~\rmd \Omega; \hspace{0.25cm}
\ell(\vv) = \int\limits_{\Gamma_t}\vv \cdot \overline{\mathbf{t}}~\rmd \Gamma
\end{equation}
where $H^1(\Omega)$ is a Hilbert space with weak derivatives up to order $1$. The weak formulation of the static problem is then given by:
\begin{equation}
\textup{find} ~\hspace{0.3cm} \uu^h \in \mathcal{U}^h \hspace{0.3cm} \textup{such~~that} \hspace{0.3cm} \forall \vv^h \in \mathcal{V}^h \hspace{0.3cm} a(\uu^h,\vv^h) = \ell(\vv^h)
\label{eqn:variweakform1}
\end{equation}
where $\mathcal{U}^h \subset \mathcal{U}$ and $\mathcal{V}^h \subset \mathcal{V}$.
\begin{figure}[htpb]
\centering
\input{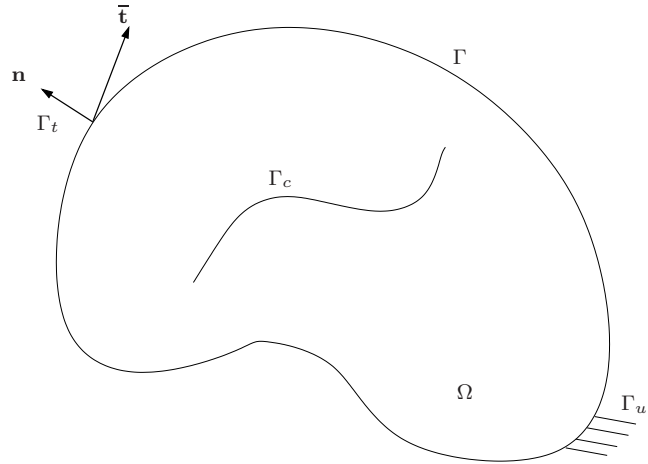}
\caption{Two-dimensional elastic body with a crack.}
\label{fig:bodywcrack}
\end{figure}

\subsection{eXtended Finite Element Method}

\begin{figure}[htpb]
\centering
\scalebox{0.7}{\input{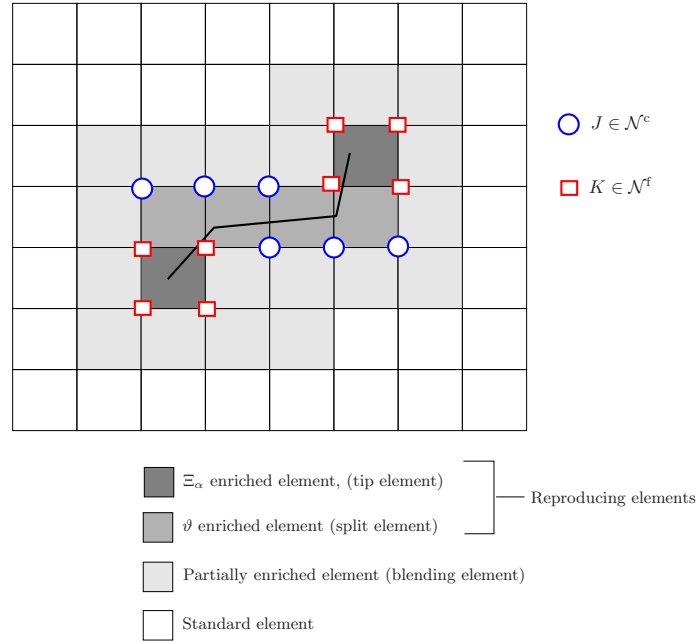}}
\caption{A typical FE mesh with an arbitrary crack. `\textit{Circled}' nodes are enriched with $\vartheta$ and `\textit{squared}' nodes are enriched with near-tip asymptotic fields.`\textit{Reproducing elements}' are the elements whose all the nodes are enriched.}
\label{fig:xfemelementcategory}
\end{figure} 

In the extended finite element method (XFEM), the classical FEM polynomial approximation space is extended by a set of functions, called the `\emph{enrichment functions}'. This is achieved through the framework of the partition of unity. By incorporating additional functions, the enriched FE basis can model singularities or internal boundaries~\cite{belytschkoblack1999}. 

\paragraph{\small Displacement approximation}
The displacement approximation can be decomposed into the standard part $\uu^h_{\rm std}$ and into an enriched part $\uu^h_{\rm xfem}$ as:

\begin{equation}
\uu^h(\xx) = \uu^h{\rm_{std}}(\xx)+ \uu^h_{\rm xfem}(\xx) 
\label{eqn:dispApproximation}
\end{equation}

In the conventional XFEM, the nonsmooth behaviour around the crack tip is represented by enhancing the conventional finite element space by appropriate functions. For the case of linear elastic fracture mechanics, two sets of functions are used: a Heaviside jump function to capture the jump across the crack faces and a set of asymptotic branch functions that span the two-dimensional asymptotic crack tip fields. A generic form of the displacement approximation to capture the displacement jump across the crack faces and to represent the singular field is given by:
\begin{equation}	
\uu^h(\xx) = \underbrace{ \sum\limits_{I \in \mathcal{N}^{\rm fem}} N_I(\xx) \qq_I}_{\uu^h{\rm_{std}}(\xx)} + \underbrace{ \sum\limits_{J \in \mathcal{N}^c} N_J(\xx) \vartheta(\xx) \aaa_J + \sum\limits_{K \in \mathcal{N}^f} N_K(\xx) \sum\limits_{\alpha=1}^4 \Xi_{\alpha}(\xx) \bb_K^\alpha}_{\uu^h_{\rm xfem}(\xx) }
\label{eqn:uS2}
\end{equation}
where $\mathcal{N}^{\rm c}$ is the set of nodes whose shape function support is cut by the crack interior (`circled' nodes in \fref{fig:xfemelementcategory}) and $\mathcal{N}^{\rm f}$ is the set of nodes whose shape function support is cut by the crack tip (`squared' nodes in \fref{fig:xfemelementcategory}). $\vartheta$ and $\Xi_{\alpha}$ are the enrichment functions chosen to capture the displacement jump across the crack surface and the singularity at the crack tip, $\qq_I$ are the standard degrees of freedom, $\aaa_J$ and $\bb^{\alpha}_K$ are the nodal degrees of freedom corresponding to functions $\vartheta$ and $\Xi_{\alpha}$, respectively.  For example, a heaviside function is used to capture the jump in the displacement, whilst $\Xi_{\alpha}$ is chosen to represent the near tip asymptotic fields. In the case of linear elastic fracture mechanics, these functions are given by:
\begin{equation}
\{\Xi_\alpha\}_{1 \le \alpha \le 4} (r,\theta)= \sqrt{r} \left\{ \sin \frac{\theta}{2}, \cos\frac{\theta}{2}, \sin\theta \sin\frac{\theta}{2},\sin\theta \cos\frac{\theta}{2} \right\}
\label{eqn:asymptotic}
\end{equation}

For orthotropic materials, the asymptotic functions given by \Eref{eqn:asymptotic} have to be modified because the material property is a function of material orientation. Asadpoure and Mohammadi~\cite{asadpouremohammadi2007,asadpouremohammadi2007a} proposed special near-tip functions for orthotropic materials as:
\begin{equation}
\{\Xi_\alpha\}_{1 \leq \alpha \leq 4} (r,\theta)= \sqrt{r} \left\{
\cos \frac{\theta_1}{2} \sqrt{g_1(\theta)},
\cos \frac{\theta_2}{2} \sqrt{g_2(\theta)},
\sin \frac{\theta_1}{2} \sqrt{g_1(\theta)},
\sin \frac{\theta_2}{2} \sqrt{g_2(\theta)} \right\}
\label{eq:orthoasymptotic}
\end{equation}
 where $(r,\theta)$ are the crack tip polar coordinates. The functions $g_i (i=1,2)$ and $\theta_i (i=1,2)$ are given by:
\begin{eqnarray}
g_j(\theta) = \left( \cos^2 \theta + \frac{\sin ^2 \theta}{e_j^2} \right), \hspace{0.5cm} j=1,2 \nonumber \\
\theta_j = \arctan \left( \frac{\tan \theta}{e_j} \right). \hspace{0.5cm} j=1,2
\end{eqnarray}
where $e_j (j=1,2)$ are related to material constants, which depend on the orientation of the material~\cite{asadpouremohammadi2007}. The local enrichment strategy introduces the following four types of elements apart from the standard elements:

\begin{itemize}
\item \textit{Split elements} are elements completely cut by the crack. Their nodes are enriched with the discontinuous function $\vartheta$.
\item \textit{Tip elements} either contain the tip or are within a fixed distance independent of the mesh size. All nodes belonging to a tip element are enriched with the near-tip asymptotic fields.
\item \textit{Tip-blending elements} are elements neighbouring tip elements. They are such that some of their nodes are enriched with the near-tip and others are not enriched at all.
\item \textit{Split-blending elements} are elements neighbouring split elements. They are such that some of their nodes are enriched with the strongly or weakly discontinuous function and others are not enriched at all.
\end{itemize}

\paragraph{\small Discretized Form}
By following Bubnov-Galerkin procedure, for arbitrary $\vv^h$, we can write the discretized form as:

\begin{equation}
\renewcommand\arraystretch{1.2}
\left[\begin{array}{*{20}c}
               \kk_{uu}  & \kk_{ua}\\
                \kk_{au}  &  \kk_{aa}
                \end{array} \right] \left\{ \begin{array} {c} \qq \\
\aaa \end{array} \right\}
= \left\{ \begin{array} {c} \ff_q \\ \ff_a \end{array} \right\}
\label{eqn:xfemKUF}
\end{equation}
where $\kk_{uu}, \kk_{aa}$ and $\kk_{ua}, \kk_{au}$ are the stiffness matrix associated with the standard FE approximation, the enriched approximation and the coupling between the standard FE approximation and the enriched approximation, respectively. 

\subparagraph{\small Remark}
The modification of the displacement approximation given by \Eref{eqn:uS2} does not introduce a new form of the discretized finite element equilibrium equations, but leads to an enlarged problem to solve.

\subsection{Difficulties in the XFEM}
Although XFEM is robust and applied to a wide variety of moving boundary problems and interface problems, the flexibility provided by this class of methods also leads to associated difficulties:

\begin{itemize}
\item{\bf Singular and discontinuous integrands} When the approximation is discontinuous or non-polynomial in an element, special care must be taken to numerically integrate over enriched elements. One potential solution for numerical integration is to partition the elements into subcells (triangles for example) aligned to the discontinuous surface in which the integrands are continuous and differentiable~\cite{belytschkoblack1999,bordasnguyen2007}. The purpose of sub-dividing into triangles is solely for the purpose of numerical integration and does not introduce new degrees of freedom. The other possible alternatives include: (a) Polar integration~\cite{b'echetminnebo2005,labordepommier2005}; (b) Complex mapping~\cite{natarajanmahapatra2009,natarajanmahapatra2010,natarajanbaiz2011a}; (c) Smoothed XFEM~\cite{bordasnatarajan2011,bordasrabczuk2010,baiznatarajan2011}; (d) Generalized quadrature~\cite{mousavisukumar2010a,mousavisukumar2010b,mousavixiao2010}; (e) Equivalent polynomials~\cite{ventura2006} and (f) Adaptive integration techniques~\cite{strouboulisbabuvska2000,liuxiao2004,xiaokarihaloo2006}.

\item{\bf Blending the different partitions of unity} The local enrichment used in the conventional element leads to oscillations in the results over the elements that are partially enriched. This pathological behaviour has been studied in great detail in the literature. Some of the proposed techniques are: Assumed strain blending elements~\cite{gracie2008,chessawang2003}, Direct coupling of the enriched and the standard regions~\cite{labordepommier2005,gracie2008,chahinelaborde2011}, Corrected or weighted XFEM and Fast integration~\cite{fries2008,venturagracie2009}, Hybrid crack element~\cite{xiaokarihaloo2003,xiaokarihaloo2007}, Hierarchical or spectral functions~\cite{tranac'onvercher2009,legaywang2005}.

\item{ \bf Ill-conditioning} The addition of enrichment functions to the FE approximation basis could result in a severely ill-conditioned stiffness matrix. The consequence of this is that the stiffness matrix can become ill-conditioned~\cite{strangfix1973,fixgulati1973,babuvskabanerjee2011,menkbordas2011}. Some of the approaches proposed in the literature include Cholesky decompsition~\cite{strangfix1973,fixgulati1973,b'echetminnebo2005} of the submatrices of the stiffness matrix, domain decomposition based preconditioner~\cite{menkbordas2011}, perturbation of the stiffness matrix~\cite{strouboulisbabuvska2000,stroubouliscopps2001} and stable GFEM/XFEM~\cite{babuvskabanerjee2011,natarajanbanerjee2012}.

\item{\bf Additional unknowns} With extrinsic enrichment, additional degrees of freedom (dofs) are introduced and the number of additional dofs depends on the number of enrichment functions and the number of such enrichments required. In case of the extrinsic enrichment, the approximation introduces additional unknowns (for example $\aaa_J$, $\bb_K^{\alpha}$, see \Eref{eqn:uS2}). These additional unknowns may increase the computational effort. By employing MLS technique on overlapping subdomains with ramp function~\cite{friesbelytschko2007} or by using crack tip elements developed for phantom node method~\cite{songareias2006,rabczukzi2008}, the total number of additional unkowns can be reduced. 
\end{itemize}

It can be seen from the literature that different approaches have been proposed to overcome the aforementioned difficulties. In addition, R\'ethor\'e \textit{et al.,}~\cite{r'ethor'eroux2010} proposed to use analytical form of the displacement field in the vicinity of the crack tip, whilst retaining the conventional form in the rest of the domain. R\'ethor\'e \textit{et al.,} used two overlapping domains to accomplish the task and the coupling between the two overlapping domains is established by matching the energy. The other approaches involve the work of Chahine and co-workers, viz., Cut-off XFEM~\cite{chahinelaborde2008} and spider XFEM~\cite{chahinelaborde2008a}. Xiao and Karihaloo~\cite{xiaokarihaloo2007} employed the hybrid crack element with XFEM to improve the accuracy of the XFEM. The crack tip region was replaced with hybrid crack element and the XFEM was used to model the crack faces with jump functions. All of the above approaches, including the conventional XFEM, requires a priori knowledge of the asymptotic fields in an analytical form. Recently, Hattori \textit{et al.,} ~\cite{hattoridiaz2012} derived new set of enrichment functions for anisotropic materials. To circumvent this problem, Menk and Bordas~\cite{menkbordas2011} proposed a numerical procedure to compute these enrichment functions, especially for anisotropic polycrystals.

\subparagraph{\small Remark}
Aforementioned difficulties are an artifact due to the asymptotic enrichments. In addition, the requirement for a priori knowledge of the asymptotic fields renders the application of the XFEM directly to heterogeneous materials for which the asymptotic fields do not exist in closed form or are very complex. Moreover, in some approaches, the computation of the stress intensity factors (SIFs) require further post-processing. Also, the introduction of the asymptotic fields, increases the complexity by requiring sophisticated numerical integration scheme to integrate over the enriched elements.

\section{Basics of the scaled boundary finite element method}\label{sbfem}
The scaled boundary finite element method (SBFEM) is a novel computational method developed by Wolf and Song~\cite{wolfsong2001}, which reduces the governing partial differential equations to a set of ordinary differential equations. The SBFEM is semi-analytical and is suitable for solving linear elliptic, parabolic and hyperbolic partial differential equations. The SBFEM relies on transforming the conventional coordinate system to the scaled boundary coordinate system. Numerical solutions are sought around the circumferential direction using conventional finite element method, whilst in the radial direction smooth analytical solutions are obtained. Like the FEM, no fundamental solution is required and like the BEM, the spatial dimension is reduced by one, since only the boundary need to be discretized, resulting in a decrease in the total degrees of freedom. 

In the SBFEM, a scaling centre $O$ is selected at a point from which the whole boundary of the domain is visible (scaling requirement). This requirement can always be satisfied by sub-structuring, i.e. dividing the structure into smaller subdomains. When modeling a cracked structure, a subdomain surrounding the crack tip is selected and the scaling centre is placed at the crack tip. The boundary of the subdomain is divided in to line elements (see \fref{fig:crkPolySBFEM}). 

\begin{figure}[htpb]
\centering
\scalebox{0.7}{\input{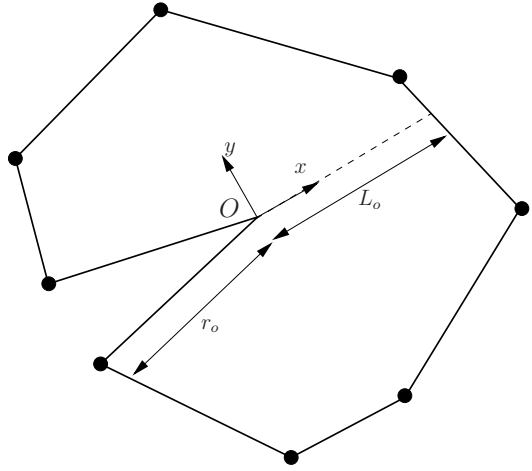}}
\caption{A cracked domain modelled by SBFEM and the definition of local coordinate system, where the `black' dots represent the nodes.}
\label{fig:crkPolySBFEM}
\end{figure}

The nodal coordinates on boundary are denoted as $\mathbf{x}_b$. As in a standard 1D isoparametric finite element,  the geometry of the element described by the coordinates $\mathbf{x}_b(\eta)$, is expressed as
\begin{equation}
 \mathbf{x}_b(\eta)=\bn(\eta)\mathbf{x}_b \label{eq:lineEleGeom}
\end{equation}
where  $\bn(\eta)$ are the shape functions. Without loss of generality, the origin of the Cartesian coordinate system is chosen at the scaling centre. The geometry of the subdomain, described by $\mathbf{x}$, is formed by  scaling the boundary (Equation~\ref{eq:lineEleGeom})
\begin{equation}
\mathbf{x} =  \xi \mathbf{x}_b(\eta) 
\label{eqn:sbccartesian}
\end{equation}
where  $\xi$ is the normalized radial coordinate running from the scaling center towards the boundary, with $\xi=$ 0 at the scaling center and $\xi=$ 1 on the boundary.  The coordinates $\xi$ and $\eta$ are the so-called scaled boundary coordinates. They  are related to the polar coordinates $r$ and $\theta$. The transformation is expressed as 
\begin{align}\label{eqn:polar}
r(\xi,\eta) & =\xi r_b(\eta) \nonumber \\
\theta(\eta)&=\arctan\dfrac{y(\eta)}{x(\eta)}
 \end{align}
where $r_b$ is the distance from the scaling centre to a point on the boundary. The transformation between the  Cartesian coordinates and the scaled boundary coordinates is similar to the coordinate transformation in constructing isoparametric finite elements. 

\subsection{Displacement approximation}
The displacements at any point is approximated by:
\begin{equation}
\uu(\xi,\eta) = \bn(\eta)\uu(\xi)
\label{eqn:dispapprox}
\end{equation}
where $\bn(\eta)$ are the shape functions of elements on the boundary and $\uu(\xi)$ is the displacement along the radial lines, represented by a set of $N$ analytical functions. By substituting \Eref{eqn:dispapprox} in the definition of strain-displacement relations, the strains $\bveps(\xi,\eta)$ are expressed as:
\begin{equation}
\bveps(\xi,\eta) = \mathbf{L}\uu(\xi,\eta)
\label{eqn:sbfemstrain}
\end{equation}
where $\mathbf{L}$ is a linear operator matrix formulated in the scaled boundary coordinates
\begin{equation}
\mathbf{L} = \bb_1(\eta) \frac{\partial}{\partial \xi} + \frac{1}{\xi}\bb_2(\eta)
\label{eqn:Loperator}
\end{equation}
and
\begin{align}
\bb_1(\eta) &= \frac{1}{det(J)} \left[ \begin{array}{cc} y_{b}(\eta)_{,\eta} & 0 \\ 0 & -x_{b}(\eta)_{,\eta} \\ - x_{b}(\eta)_{,\eta} & y_{b}(\eta)_{,\eta} \end{array} \right], \nonumber \\
\bb_2(\eta) &= \frac{1}{det(J)} \left[ \begin{array}{cc} -y_b & 0 \\ 0 & x_b \\ x_b & y_b \end{array} \right],
\end{align}
The determinant of the Jacobian matrix is: 
\begin{equation}det(J) = x_b(\eta)y_b(\eta)_{,\eta} - y_b(\eta)x_b(\eta)_{,\eta}.
\end{equation}
where $x_b(\eta)$ and $y_b(\eta)$ are given by \Eref{eq:lineEleGeom}. The stresses $\bvsig(\xi,\eta)$ are then given by:
\begin{equation}
\bvsig(\xi,\eta) = \dd \bigb_1(\eta)\uu(\xi)_{,\xi} + \frac{1}{\xi} \dd \bigb_2(\eta) \uu(\xi)
\label{eqn:sbfemstress}
\end{equation}
where in the above equation, the definition of strain and the linear operator matrix given by \Eref{eqn:sbfemstrain} - (\ref{eqn:Loperator}) are used with
\begin{align}
\bigb_1(\eta) & = \bb_1(\eta)\bn(\eta) \nonumber \\
\bigb_2(\eta)	& = \bb_2(\eta)\bn(\eta)_{,\eta}
\end{align}
Upon substituting the \Erefs{eqn:sbfemstress} and (\ref{eqn:sbfemstrain}) in the virtual work statement and following the derivation in~\cite{wolfsong2001,deekswolf2002}, the virtual work yields
\begin{equation}
\begin{split}
\delta \uu_b^{\rm T} \left((\ee_o \xi \uu(\xi)_{,\xi} + \ee_1^{\rm T} \uu(\xi))|_{\xi=1}  -  \ff\right) \\
- \int\limits_0^1 \delta \uu(\xi)^{\rm T} \left( \ee_o \xi^2 \uu(\xi)_{,\xi\xi} + (\ee_o + \ee_1^{\mathrm{T}} - \ee_1)\xi\uu(\xi)_{,\xi} - \ee_2\uu(\xi) \right)~\rmd \xi = 0
\end{split}
\label{eqn:sbfemvirtual}
\end{equation}
where $\ff$ is the equivalent boundary nodal forces and $\uu_b$ is the nodal displacement vector. By considering the arbitrariness of $\delta \uu(\xi)$, the following ODE is obtained:
\begin{equation}
\ee_o \xi^2 \uu(\xi)_{,\xi\xi} + (\ee_o + \ee_1^{\mathrm{T}} - \ee_1)\xi\uu(\xi)_{,\xi} - \ee_2\uu(\xi) = 0
\label{eqn:governODEsbfem}
\end{equation}
where $\ee_o, \ee_1$ and $\ee_2$ are known as the coefficient matrices and are given by:
\begin{align}
\ee_o &= \int_\eta \bigb_1(\eta)^{\rm T} \dd \bigb_1(\eta) det(J)~\rmd\eta, \nonumber \\
\ee_1 &= \int_\eta \bigb_2(\eta)^{\rm T} \dd \bigb_1(\eta) det(J)~\rmd\eta, \nonumber \\
\ee_2 &= \int_\eta \bigb_2(\eta)^{\rm T} \dd \bigb_2(\eta) det(J)~\rmd\eta.
\label{eqn:coeffmat}
\end{align}
The boundary nodal forces are related to the displacement functions by
\begin{equation}
\ff^{^{\rm sbfem}}  =( \ee_o \xi \uu(\xi)_{,\xi} + \ee_1^{\rm T} \uu(\xi))|_{\xi=1} \label{eqn:nodalforce}
\end{equation}

\subparagraph{\small Remark}
For a 2 noded line element, the coefficient matrices given by \Eref{eqn:coeffmat} can be written explicitly. This is given in Appendix A. A MATLAB routine (EleCoeffMatrices2NodeEle.m) is given in Appendix B to compute the coefficient matrices.

\subsection{Computation of the stiffness matrix}
\Eref{eqn:governODEsbfem} is a homogeneous second-order ordinary differential equation, the general solution is given by~\cite{wolfsong2001,deekswolf2002}:
\begin{equation}
\uu(\xi) = \pphi  \xi^{-\boldsymbol{\lambda}} \mathbf{c} = \sum_i c_i \xi^{-\lambda_i} \phi_i
\label{eqn:dispgeneralsoln}
\end{equation}
where $\pphi_i$ and $\lambda_i$  can be interpreted as the deformation modes and the corresponding scaling factors that closely satisfy internal equilibrium in the $\xi$ direction~\cite{wolfsong2001,deekswolf2002}, and the constants $c_i$ are the contribution of each mode. The deformation modes are obtained from the following quadratic eigenvalue problem obtained by substituting \Eref{eqn:dispgeneralsoln} into \Eref{eqn:governODEsbfem}  as:

\begin{equation}
\left( \lambda_i^2  \ee_o  - \lambda_i (\ee_1^{\mathrm{T}} - \ee_1) - \ee_2\right) \phi_i = 0
\label{eqn:lameqn}
\end{equation}
This equation yields $2N$ modes with pairs of eigenvalues ($\lambda$, -$\lambda$). Only the $N$ modes with negative real parts of eigenvalue, which leads to finite displacements at the scaling centre, are selected. On the boundary, the displacement $\mathbf{u}_b = \mathbf{u} (\xi =1)$ is given by (\Eref{eqn:dispgeneralsoln}):
\begin{equation}
\mathbf{u}_b = \pphi \mathbf{c}
\end{equation}
which permits the integration constants $\mathbf{c}$ to be determined from displacement on the boundary
\begin{equation}
\mathbf{c} = \pphi^{-1} \uu_b \label{eqn:intgconst}
\end{equation}
The nodal force on the boundary is obtained by substituting \Eref{eqn:dispgeneralsoln} into \Eref{eqn:nodalforce} as:
\begin{equation}
\ff^{^{\rm sbfem}} = \left(- \ee_o  \pphi \boldsymbol{\lambda} + \ee_1^{\rm T} \pphi \right)\mathbf{c}
\label{eqn:sbfemforce1}
\end{equation}
Substituting \Eref{eqn:intgconst} into \Eref{eqn:sbfemforce1}  yields:
\begin{equation}
\kk^{^{\rm sbfem}} \uu_b = \ff^{^{\rm sbfem}} 
\label{eqn:discretizesbfem}
\end{equation}
where $\kk^{^{\rm sbfem}} $ is the stiffness matrix given by:
\begin{equation}
\kk^{^{\rm sbfem}}  = - \ee_o \pphi\boldsymbol{\lambda}\pphi^{-1} + \ee_1^{\rm T}
\label{eqn:sbfemkmat}
\end{equation}

Alternatively, the quadratic eigen-problem in \Eref{eqn:lameqn} is rewritten as a standard eigen-problem of the matrix 
\begin{equation}
\mathbf{Z} = \left[ \begin{array}{cc} \ee_o^{-1} \ee_1^{\rm T} & -\ee_o^{-1} \\
-\ee_2 + \ee_1 \ee_o^{-1} \ee_1^{\rm T} & -\ee_1  \ee_o^{-1} 
\end{array}\right]
\label{eq:Zmatrix}
\end{equation}
The eigenvalues and eigenvectors are partitioned as
\begin{equation}
\mathbf{Z}  
\left[ \begin{array}{cc} \mathbf{v}_{11}  & \mathbf{v}_{12} \\
\mathbf{v}_{21} & \mathbf{v}_{22}
\end{array}\right]
=
\left[ \begin{array}{cc} \mathbf{v}_{11}  & \mathbf{v}_{12} \\
\mathbf{v}_{21} & \mathbf{v}_{22}
\end{array}\right]
\left[ \begin{array}{cc} \boldsymbol{\lambda}  &0 \\
0 & -\boldsymbol{\lambda}
\end{array}\right]
\end{equation}
where $\boldsymbol{\lambda}$ consists of all the eigenvalues with negative real parts. The displacement  solution is expressed as (\Eref{eqn:dispgeneralsoln}) 
\begin{equation}
\uu(\xi) =  \mathbf{v}_{11}  \xi^{-\boldsymbol{\lambda}} \mathbf{c} 
\end{equation}
The boundary nodal force (\Eref{eqn:sbfemforce1}) is written as
\begin{equation}
\ff^{^{\rm sbfem}} = \mathbf{v}_{21} \mathbf{c}
\end{equation}
The stiffness matrix is obtained as
\begin{equation}
\kk^{^{\rm sbfem}}  = \mathbf{v}_{21}  \mathbf{v}_{11}^{-1} 
\label{eqn:sbfemkmat-b}
\end{equation}

\subparagraph{\small Remark}
A MATLAB routine (getSBFEMStiffMat.m) is given in Appendix B to compute the stiffness matrix, the deformation modes and the corresponding scaling factors.  It is based on Equations~(\ref{eq:Zmatrix}) - (\ref{eqn:sbfemkmat-b}).


\section{eXtended Scaled Boundary Finite Element Method}
\label{xsbfemcou}
The main idea of the proposed xSBFEM is to approximate the nonsmooth behaviour around the crack tip by a semi-analytical approach described in the previous section. In the proposed approach, the domain of $\Omega$ is partitioned into three non-overlapping regions $(\Omega = \Omega^{^{\rm fem}} \cup \Omega^{^{\rm xfem}} \cup \Omega^{^{\rm sbfem}})$ (see \fref{fig:xfemsbfemcoup}):

\begin{itemize}
\item $\Omega^{^{\rm fem}}$ - contains a set of elements that are not intersected by the discontinuity surface. In this region the conventional FE approximation is employed to represent the smooth behaviour. 
\item $\Omega^{^{\rm xfem}}$ - contains a set of elements that are completely cut by the crack. Also called as the split element in XFEM terminology. The support of their nodes are enriched with a Heaviside function.
\item $\Omega^{^{\rm sbfem}}$ - set of elements that define the SBFEM domain. The nonsmooth behaviour around the crack tip is modelled by a semi-analytical approach.
\end{itemize}


\begin{figure}[htpb]
\centering
\scalebox{0.5}{\input{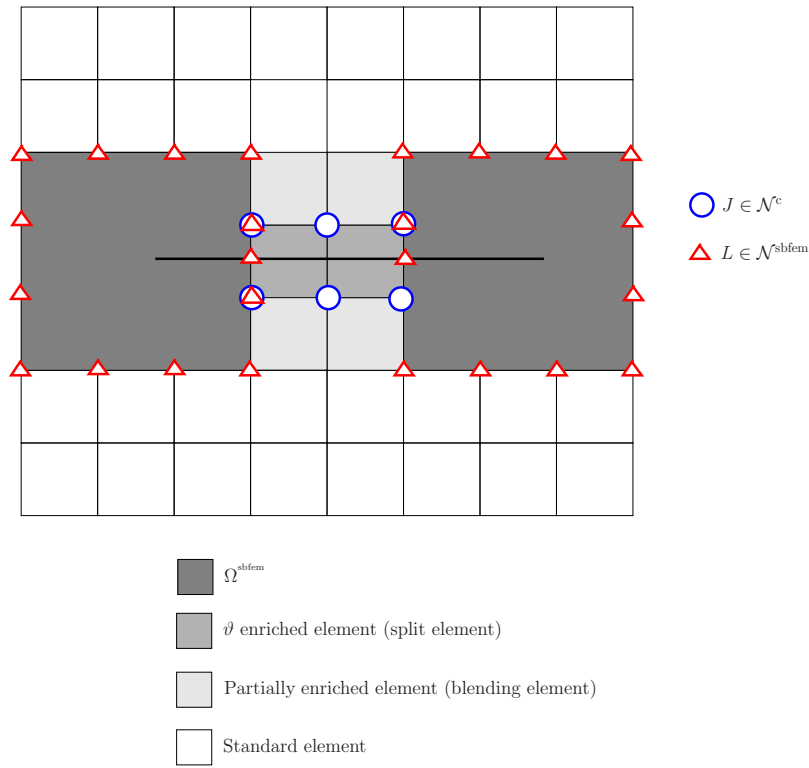}}
\caption{A typical FE mesh with a crack and different regions identified to represent the behaviour of the problem. In this descriptive figure, one layer of elements around the crack tip is replaced with the SBFEM domain. `\textit{Circled}' nodes are enriched with $\vartheta$ and `\textit{triangled}' nodes are the nodes on the boundary of the $\Omega^{^{\rm sbfem}}$. Note that at the crack mouth, SBFEM requires two nodes, one on either side of the crack, whereas the discontinuity is captured by the enrichment function $\vartheta$ in case of the XFEM.}
\label{fig:xfemsbfemcoup}
\end{figure}

\subparagraph{\small Remark}
The SBFEM domain is constructed from the underlying FE mesh. After identifying a region around the crack tip, the information on the boundary of this region alone is required. The degrees of freedom of the unused nodes in the FE mesh are taken care of during the solution process.

\subsection{Salient features of the proposed approach} The proposed approach aims at improving the accuracy of the XFEM in modeling problems with strong discontinuity. It shares the idea of representing the nonsmooth behaviour in the vicinity of the crack tip by analytical form similar to the work of R\'ethor\'e \textit{et al.,}~\cite{r'ethor'eroux2010} and Chahine \textit{et al.,}~\cite{chahinelaborde2008} but differs in basic approach and implementation. Some of the salient features are listed below:

\begin{itemize}
\item An underlying FE mesh is used to define the SBFEM domain. 
\item No a priori knowledge of the asymptotic fields is required.
\item In the vicinity of the crack tip, the integration is performed along the boundary of the SBFEM domain.
\item No special technique is required to match the different partitions of unity.
\item The existing approach can easily be incorporated into any existing FEM/XFEM code.
\end{itemize}

\subsection{Displacement approximation and the stiffness matrix}
The displacement approximation can be decomposed into the standard part $\uu^h_{\rm std}$ and into an enriched part $\uu^h_{\rm xfem}$ as given by \Eref{eqn:uS2} without asymptotic enrichments. In the proposed approach, the standard part is further decomposed into the standard FEM $\uu^h_{\rm fem}$ and into a semi-analytical part $\uu^h_{\rm sbfem}$. The displacement approximation is then given by:

\begin{equation}
\renewcommand{\arraystretch}{2}
\uu^h(\xx) = \left\{ \begin{array}{lr} \sum\limits_{I \in \mathcal{N}^{\rm fem}} N_I(\xx) \qq_I + \sum\limits_{J \in \mathcal{N}^c} N_J(\xx) H(\xx) \aaa_J & \xx \in \Omega^{^{\rm fem}} \cup \Omega^{^{\rm xfem}} \\
 \sum\limits_{K \in \mathcal{N}^{\rm sbfem}} N_K(\eta)  \uu_{bK}(\xi) & \xx = f(\xi,\eta) \in \Omega^{^{\rm sbfem}} \end{array} \right.
 \label{eqn:dispApprox}
\end{equation}
where $\mathcal{N}^{\rm fem}$ is the set of all nodes in the FE mesh that belongs to $\Omega^{^{\rm fem}}$, $\mathcal{N}^{\rm sbfem}$ is the set of nodes in the FE mesh that belongs to $\Omega^{^{\rm sbfem}}$ and $\mathcal{N}^{\rm xfem}$ is the set of nodes that are enriched with the Heaviside function $H$. $N_I$ and $N_J$ are the standard finite element functions, $\qq_I$ and $\aaa_J$ are the standard and the enriched nodal variables associated with node $I$ and node $J$, respectively. $\uu_{bK}$ are the nodal variables associated with the nodes on the boundary of $\Omega^{^{\rm sbfem}}$ and $N_K$ are the standard FE shape functions defined on the scaled boundary coordinates in $\Omega^{^{\rm sbfem}}$ (see Section \ref{sbfem}). 

For $\Omega = \Omega^{^{\rm fem}} \cup \Omega^{^{\rm xfem}}$, the discretized is given by \Eref{eqn:xfemKUF}. For $\Omega = \Omega^{^{\rm sbfem}}$, the discretized form is given by \Eref{eqn:discretizesbfem} and the stiffness matrix is given by \Eref{eqn:sbfemkmat}. Note that $\uu_b$ in \Eref{eqn:discretizesbfem} are the nodes on the boundary of SBFEM domain $\Omega^{^{\rm sbfem}}$ that share with the boundary of $\Omega^{^{\rm fem}}$ and $\Omega^{^{\rm xfem}}$. These nodes are identified from the underlying FE mesh. Although, the examples shown here are for structured mesh, the proposed technique is not limited to a structured mesh. Next, we describe the coupling of the three regions. 


\subsection{Blending different partitions of unity} 
 Before discussing the coupling of different regions, let us define the notations used for the displacement variables in different regions.
 
 \begin{itemize}
 \item $\qq_I$ - standard degrees of freedom in $\Omega^{^{\rm fem}}$ and $\uu_{qS}$ are the standard degrees of freedom in $\Omega^{^{\rm fem}} \cap \Omega^{^{\rm sbfem}}$.
 \item $\aaa_I$ - enriched degrees of freedom in $\Omega^{^{\rm xfem}}$ corresponding to the enrichment function $\vartheta$ that describes the jump in the displacement across the discontinuity surface.
 \item $\uu_{bK}$ - degrees of freedom in $\Omega^{^{\rm sbfem}}$. This is further decomposed into two parts: $\uu_{bK}^1$ are the degrees of freedom for the nodes in $\Omega^{^{\rm fem}} \cap \Omega^{^{\rm sbfem}}$ and $\uu_{bK}^2$ are the degrees of freedom for the nodes in $\Omega^{^{\rm sbfem}} \cap \Omega^{^{\rm xfem}}$.
 \item $\uu_{xF}$ - contains the standard and the enriched degrees of freedom for the nodes in $\Omega^{^{\rm xfem}}$ that share with the $\Omega^{^{\rm sbfem}}$.
 \end{itemize}

\paragraph{ \small  Coupling of~$\Omega^{^{\rm fem}}$, $\Omega^{^{\rm xfem}}$ and $\Omega^{^{\rm sbfem}}$} For the nodes that on the boundary of $\Omega^{^{\rm fem}}$ and $\Omega^{^{\rm sbfem}}$, no special coupling technique is required. As the nodal unknown coefficients are required to be continuous across the boundary, the unknown coefficients from the SBFEM and the FEM are the same and are assembled in the usual way. A similar procedure is followed when assembling the $\Omega^{^{\rm fem}}$ and $\Omega^{^{\rm xfem}}$ to the global stiffness and to the global force vector. 

\paragraph{\small  Coupling of~$\Omega^{^{\rm sbfem}}$ and $\Omega^{^{\rm xfem}}$}
In case of the SBFEM, a crack in a body is modelled by an open domain as shown in \fref{fig:crkPolySBFEM} with two nodes, one on either side of the crack, whilst, by using a heaviside function, a crack is modeled in the XFEM. The addition of a heaviside function introduces additional degrees of freedom that accounts for the jump in the displacement across the crack surface. \fref{fig:xfemsbfemcouponelement} shows a typical approach to model crack using the XFEM and the SBFEM. The displacement approximation for the XFEM and the SBFEM is given by:
\begin{eqnarray}
\uu^h_{\rm xfem} &=& \sum\limits_{I=1}^4 N_I(x,y) \qq_I + \sum\limits_{I=1}^4 N_I(x,y) \left[ H(\xx) - H(\xx_I) \right] \aaa_I \nonumber \\
\uu^h_{\rm sbfem} &=& \sum\limits_{I=1}^{n_b} N_I(\eta) \uu_{bI}(\xi)
\label{eqn:enrichApprox}
\end{eqnarray}
where $n_b$ is the total number of nodes on the boundary of $\Omega^{^{\rm sbfem}}$. As the displacements are required to be continuous, the displacements from the XFEM domain are matched to the displacements from the SBFEM domain by using a transformation matrix (see below). The continuity of the displacements are enforced only on the side that is common to both the $\Omega^{^{\rm xfem}}$ and $\Omega^{^{\rm sbfem}}$.

\paragraph{\small Numerical illustration}
The coupling between these two regions is illustrated below for a simple case where the crack is straight and the $\Omega^{^{\rm sbfem}}$ is to the right of $\Omega^{^{\rm xfem}}$ (see \fref{fig:xfemsbfemcouponelement}). The displacement approximation for the XFEM and the SBFEM is given by \Eref{eqn:enrichApprox}. As explained earlier, to ensure compatibility of displacements, the displacements from the SBFEM domain are transformed to the displacements in the XFEM domain. This is done on the side that is common to both the regions, viz., sides connecting 2-3 in the XFEM domain and B-A: F-E in the SBFEM domain (see \fref{fig:xfemsbfemcouponelement}). From the XFEM approximation for the displacements, the displacements at the mouth of the crack, ie., at points $`F'$ and $`A'$ can be written as:

\begin{eqnarray}
\uu_F(x,y) &=& \sum\limits_{I=2,3} N_I(\xx_F) \qq_I + \sum\limits_{I=2,3} N_I(\xx_F) \left[H(\xx_F) - H(\xx_I) \right] \aaa_I \nonumber \\
&=& \sum\limits_{I=2,3} N_I(\xx_F) \qq_I +  2 N_2(\xx_F)  \aaa_2 \nonumber \\
\uu_A(x,y) &=& \sum\limits_{I=2,3} N_I(\xx_A) \qq_I + \sum\limits_{I=2,3} N_I(\xx_A) \left[H(\xx_A) - H(\xx_I) \right] \aaa_I \nonumber \\
&=& \sum\limits_{I=2,3} N_I(\xx_A) \qq_I -  2 N_3(\xx_A)  \aaa_3
\label{eqn:xfemsbfemdisp}
\end{eqnarray}

\begin{figure}[htpb]
\centering
\scalebox{0.85}{\input{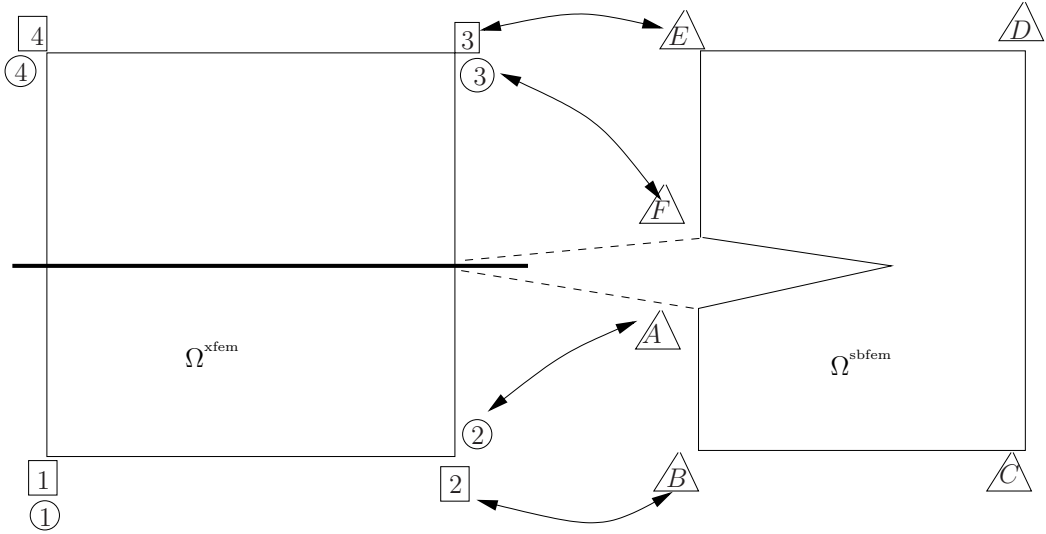}}
\caption{Coupling the XFEM and the SBFEM domain. \textit{`circled'} nodes are the enriched nodes. `\textit{Dashed}' line represent the crack.}
\label{fig:xfemsbfemcouponelement}
\end{figure}

Note that in \Eref{eqn:xfemsbfemdisp}, the displacements $\uu_F$ and $\uu_A$ are expressed in terms of nodal displacements $\qq_I$ and the enriched displacements $\aaa_I$. Also note that these displacement should match with the displacements from the SBFEM domain. A transformation matrix can be derived from the above equation that translates the SBFEM nodal unknown displacements (viz, $\uu_A, \uu_B, \uu_E, \uu_F$) to the XFEM nodal unknown displacements (viz, $\qq_2, \qq_3, \aaa_2, \aaa_3$). In matrix form,

\begin{equation}
\left\{ \begin{array}{c} \uu_B \\ \uu_E \\ \uu_A \\ \uu_F \end{array} \right\} = \left[ \begin{array}{cccc} \mathbf{I} & 0 & 0 & 0 \\ 0 & \mathbf{I} & 0 & 0  \\ N_2(\xx_A) & N_3(\xx_A) & 0 & -2N_3(\xx_A) \\ N_2(\xx_F) & N_3(\xx_F) & 2N_2(\xx_F) & 0 \end{array} \right] \left\{ \begin{array}{c} \qq_2 \\ \qq_3 \\ \aaa_2 \\ \aaa_3 \end{array} \right\} \\
\implies \uu^2_{bK} = \mathbf{T} \uu^h_{xF}
\label{eqn:xfemsbfemT}
\end{equation} 
where $\mathbf{I}$ is the identify matrix and $\mathbf{T}$ is the transformation matrix that relates the displacements on the side that is common to $\Omega^{^{\rm sbfem}}$ and $\Omega^{^{\rm xfem}}$. To ensure the compatibility of the displacements and to assemble the stiffness matrix to the global stiffness matrix, the stiffness matrix in $\Omega^{^{\rm sbfem}}$ and the force vector is rearranged as:
\begin{equation}
\left[ \begin{array}{cc} \kk_{aa} & \kk_{ab} \\ \kk_{ba} & \kk_{bb} \end{array} \right] \left\{ \begin{array}{c} \uu_{bK}^1 \\ \uu_{bK}^2 \end{array} \right\} = \left\{ \begin{array}{c} \ff_a \\ \ff_b\end{array} \right\}
\end{equation}
where $\kk_{aa}$, $\uu_{bK}^1 $ and $\ff_a$ corresponds to the stiffness matrix, the displacement vector and the force vector for the nodes that does not share with the element that is completely cut by the discontinuity surface in the $\Omega^{^{\rm xfem}}$, $\kk_{bb}$, $\uu_{bK}^2$ and $\ff_b$ are the stiffness matrix, the displacement vector and the force vector that share with the element that is completely cut by the discontinuity surface and $\kk_{ab}, \kk_{ba}$ are the coupling matrices. Now, applying the transformation matrix, we get
\begin{equation}
\left[ \begin{array}{cc} \kk_{aa} & \kk_{ab}\mathbf{T} \\ \mathbf{T}^{\rm T} \kk_{ba} & \mathbf{T}^{\rm T}\kk_{bb} \mathbf{T} \end{array} \right] \left\{ \begin{array}{c} \uu_{qS} \\ \uu_{xF} \end{array} \right\} = \left\{ \begin{array}{c} \ff_a \\ \mathbf{T}^{\rm T} \ff_b\end{array} \right\}
\end{equation}


\subparagraph{\small Remark}
The transformation matrix $\mathbf{T}$ requires only the evaluation of shape functions at the crack mouth, corresponding to the edge that is common to the $\Omega^{^{\rm xfem}}$ and $\Omega^{^{\rm sbfem}}$.

\section{Calculation of the stress intensity factors} \label{sifcalculation}


\subsection{SIF from semi-analytical solution of displacements} 
The relation between the stress intensity factors (SIFs) and the crack tip displacement fields for an isotropic material at a coordinate system as shown in \fref{fig:crkPolySBFEM}. The displacement field corresponding to the singular stress terms (superscript $s$ for singular) is given by:

\begin{eqnarray}
u^s_x(r,\theta) = \frac{K_I}{G} \sqrt{\frac{r}{2\pi}} \cos \frac{\theta}{2} \left[ \frac{1}{2}(\kappa-1) + \sin^2 \frac{\theta}{2} \right] + \frac{K_{II}}{G} \sqrt{\frac{r}{2\pi}} \sin\frac{\theta}{2}\left[ \frac{1}{2}(\kappa+1) + \cos^2\frac{\theta}{2} \right] \nonumber \\
u^s_y(r,\theta) = \frac{K_I}{G} \sqrt{\frac{r}{2\pi}} \sin \frac{\theta}{2} \left[ \frac{1}{2}(\kappa-1) + \cos^2 \frac{\theta}{2} \right] + \frac{K_{II}}{G} \sqrt{\frac{r}{2\pi}} \cos\frac{\theta}{2}\left[ \frac{1}{2}(\kappa+1) + \sin^2\frac{\theta}{2} \right] 
\label{eqn:asymsol}
\end{eqnarray}
where $K_I$ and $K_{II}$ are the  mode-I and mode-II SIFs,  $G$ is the shear modulus, $\nu$ is the Poisson's ratio and $\kappa$ is
\begin{equation}
\kappa = \left\{ \begin{array}{cc} 3-4\nu & \textup{for plate strain} \\ (3-\nu)/(1+\nu) & \textup{for plane stress} \end{array} \right.
\end{equation}
Consider the crack mouth opening displacement with $r=r_o$ and $\theta=\pm\pi$, we have
\begin{equation}
\left\{ \begin{array}{c} K_I \\ K_{II} \end{array} \right\} = \frac{2G}{\kappa+1} \sqrt{\frac{2\pi}{r_o}} \left\{ \begin{array}{c} u^s_y(r_o ,+\pi)-u^s_y(r_o,-\pi)  \\ u^s_x(r_o,+\pi)-u^s_x(r_o,-\pi) \end{array} \right\}
\label{eqn:SIFfromDisp}
\end{equation}

The singular stress terms can be identified by examining the eigenvalues of the terms in the semi-analytical solution (\Eref{eqn:dispgeneralsoln}) obtained in the scaled boundary finite element method. For the case of a crack in a homogeneous subdomain, the eigenvalues of the singular terms are approximately equal to $0.5$. The displacement functions of the singular term is expressed as
\begin{equation}
\uu^{s} (\xi) = \xi^{0.5} \sum\limits_{i=I,II} c_i  \phi_i \label{eqn:singulardisp}
\end{equation}
The corresponding nodal displacement on the boundary ($\xi = 1$) is written as
\begin{equation}
\uu_b^{s} =  \sum\limits_{i=I,II} c_i  \phi_i
\end{equation}
The values of the displacement of the singular stress term at the two nodes on the crack mouth  are extracted from $\uu_b^{s}$ and substituted into the last term in \Eref{eqn:SIFfromDisp} leading to the stress intensity factors  $K_I$ and $K_{II}$.

\subsection{SIF from semi-analytical solution of stresses}
The stress intensity factors $K_I$ and $K_{II}$ can be extracted from their definitions based on singular stresses  $\sigma^s_{yy}$ and   $\sigma^s_{xy}$ at the crack front $\theta = 0$ 
\begin{equation}
\left\{ \begin{array}{c} K_I \\ K_{II} \end{array} \right\} = \sqrt{2\pi r} \left\{ \begin{array}{c} \sigma^s_{yy} (\theta = 0) \\  \sigma^s_{xy} (\theta=0) \end{array} \right\}
\label{eqn:sifdefn}
\end{equation}

The singular stress modes are evaluated by using \Eref{eqn:sbfemstress} and the corresponding displacement modes in  \Eref{eqn:singulardisp} 
 \begin{equation}
\Psi_i(\eta) = \left\{ \begin{array}{c} \Psi_{xx}(\eta) \\ \Psi_{yy}(\eta) \\ \Psi_{xy}(\eta) \end{array} \right\}_i = \dd[ \lambda_i \bigb_1(\eta) + \bigb_2(\eta)]\phi_i 
\end{equation}

As in the stress recovery in the finite element method, the computation is performed element-by-element at the Gauss integration point of  an order lower  by one than that for full integration. This yields an array of values of the singular stress modes at discrete values of angular coordinate $\theta$ (\Eref{eqn:polar}). The values at the crack front $\theta=0$ are obtained by interpolation and denoted as   $\Psi_i(\theta=0)$. Along the crack front, $\xi = r/L_o$ (\Eref{eqn:polar}) applies and the singular stresses  are expressed as
\begin{equation}
\bvsig(r,\theta=0) = \sqrt{\dfrac{L_o}{r}}\sum\limits_{i=I, II} c_i  \Psi_i(\theta=0)
\label{eqn:sbfemstressparts}
\end{equation}

Substituting the stress components $\sigma_{yy}$ and $\sigma_{xy}$ obtained from \Eref{eqn:sbfemstressparts} into \Eref{eqn:sifdefn} yields the stress intensity factors

\begin{equation}
\left\{ \begin{array}{c} K_I \\ K_{II} \end{array} \right\} =  
\sqrt{2\pi L_o } \left\{ \begin{array}{c}  \sum\limits_{i=I,II} c_i \Psi_{yy} (\theta=0)_i \\  \sum\limits_{i=I,II} c_i \Psi_{xy}(\theta=0)_i \end{array} \right\}
\label{eqn:sifdefnmodes1}
\end{equation}

\subparagraph{\small Remark} 
The stress singularity is handled analytically and the standard stress recovery techniques are employed in extracting the stress intensity factors.

\section{Numerical Examples} \label{numexample}
In this section, we illustrate the effectiveness of the proposed method by solving a few benchmark problems taken from linear elastic fracture mechanics. We first consider an infinite plate under tension, then a plate with an edge crack loaded with far field tension and shear. Then as a last example, we consider an orthotropic plate with a central crack and an edge crack with far field tension. In this study, unless otherwise mentioned, bilinear quadrilateral elements are used. The results from the proposed method are compared with those of the standard XFEM and the corrected XFEM.

\subsection{Infinite plate under remote tension}
In the first example, a plate with a crack loaded in remote tension is considered. Along the boundary ABCD (see \fref{fig:griffithProblem}), the closed form near tip displacements are imposed, given by (for plane strain conditions, see \Eref{eqn:asymsol}):
\begin{eqnarray}
u_x(r,\theta) &=& \frac{2(1+\nu)}{\sqrt{2\pi}}\frac{K_I}{E} \sqrt{r}\cos\frac{\theta}{2} \left(2-2\nu-\cos^2\frac{\theta}{2} \right), \nonumber \\
u_y(r,\theta) &=& \frac{2(1+\nu)}{\sqrt{2\pi}}\frac{K_I}{E} \sqrt{r}\sin\frac{\theta}{2} \left(2-2\nu-\cos^2\frac{\theta}{2} \right)
\label{eqn:modeIcrk}
\end{eqnarray}
where $K_I = \sigma \sqrt{\pi a}$ denotes the stress intensity factor, $\sigma$ is the remote tension, $\nu$ is the Poisson's ratio and $E$ is the Young's modulus. All simulations are performed with $\sigma=$ 10$^4$ N/m$^2$ on a square with sides of length 10m.


\begin{figure}[htpb]
\centering
\scalebox{0.45}{\input{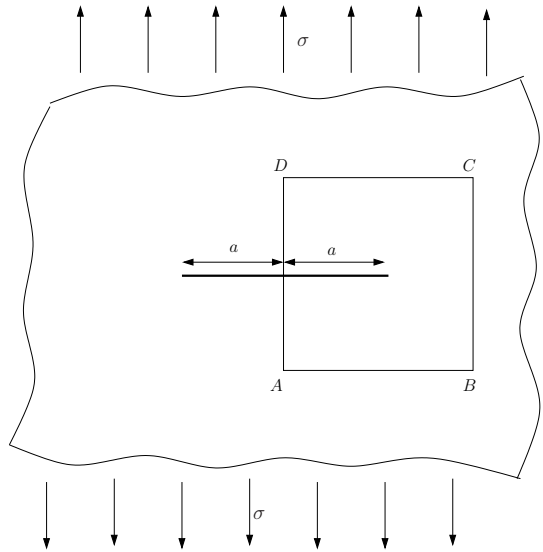}}
\caption{Infinite plate with a center crack under remote tension. Region $A-B-C-D$ is chosen as the computational domain and analytical displacements are prescribed along $A-B-C-D$ given by \Eref{eqn:modeIcrk}.}
\label{fig:griffithProblem}
\end{figure}

\begin{figure}[htpb]
\centering
\includegraphics[scale=0.65]{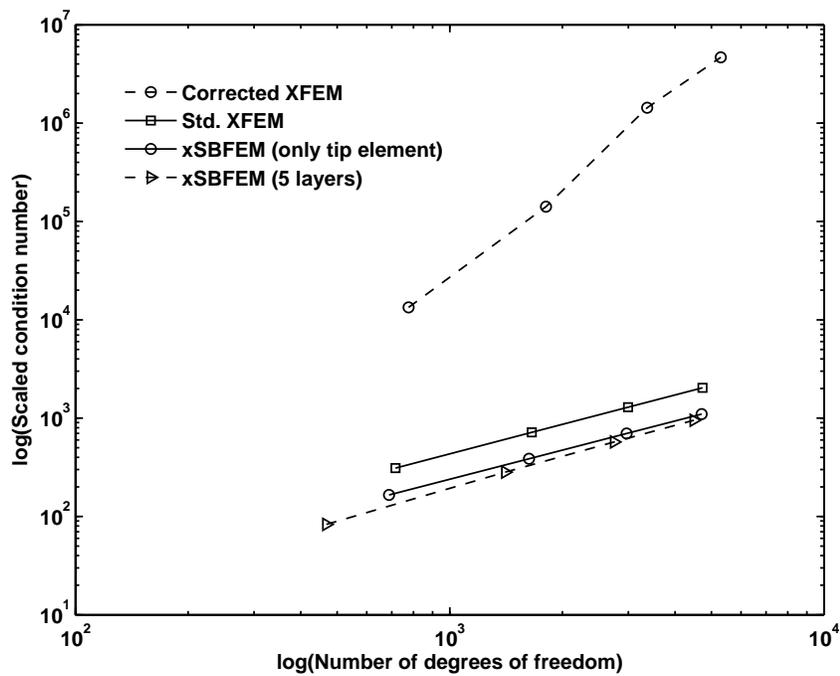}
\caption{Scaled condition number evaluated for different mesh densities. Std. XFEM refers to the conventional XFEM with `\textit{topological}' enrichment and the Corrected XFEM refers to XFEM with `\textit{geometric}' enrichment with blending correction.}
\label{fig:scaledcondk}
\end{figure}

\begin{figure}[htpb]
\centering
\includegraphics[scale=0.65]{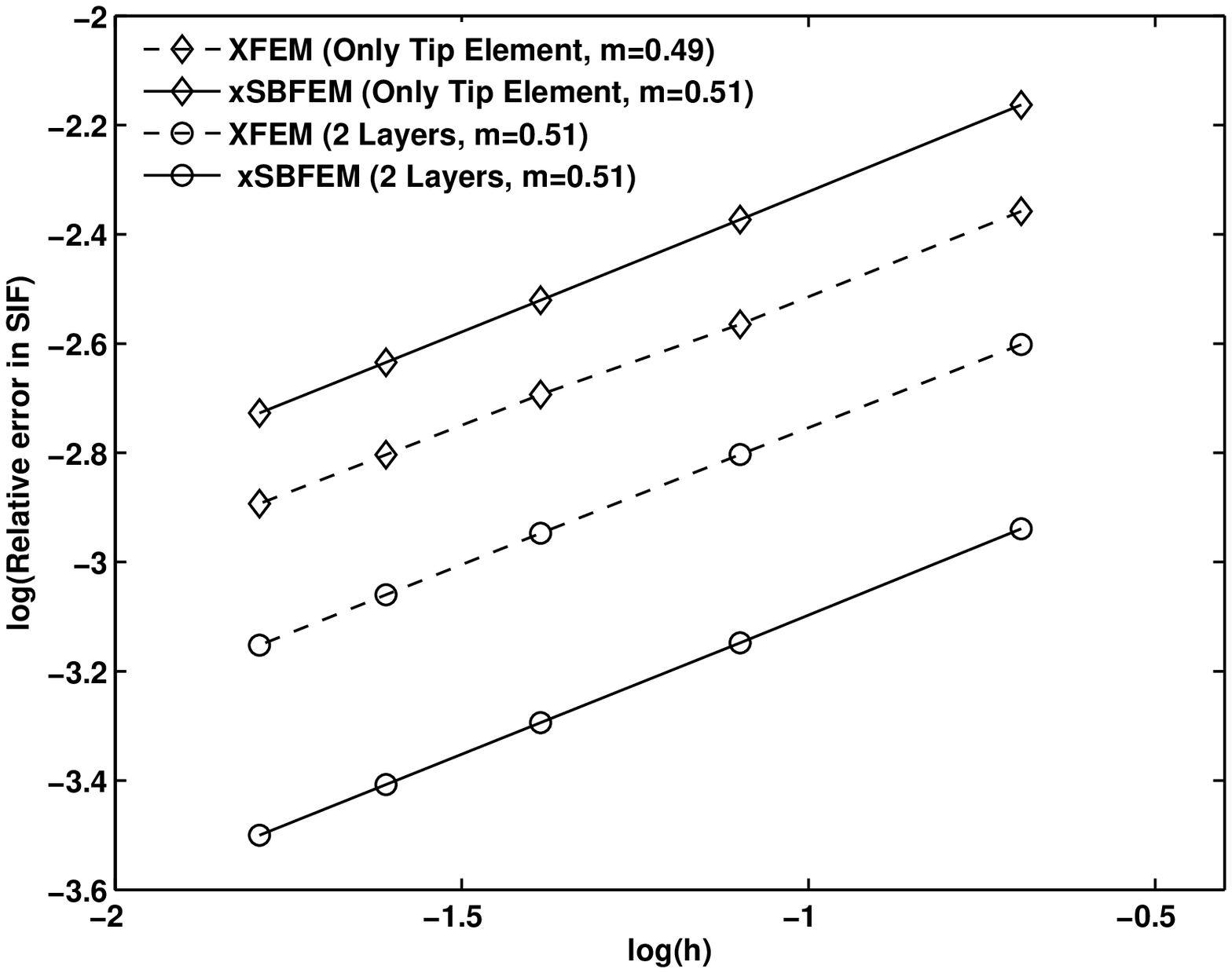}
\caption{Infinite plate under far field tension: Convergence results in the SIF. The rate of convergence is also given in the figure.}
\label{fig:conveGriffithJint}
\end{figure}

The scaled condition number $\mathfrak{K}$ of the transformed stiffness matrix $(\tilde{\mathbf{K}} = \mathbf{P}^{\rm T} \mathbf{K} \mathbf{P})$ with increasing total degrees of freedom is shown in \fref{fig:scaledcondk}. The matrix $\mathbf{P}$ is a preconditioner and it is chosen such that the condition number of the transformed matrix $\tilde{\mathbf{K}}$ is smaller than the original matrix. In this study, the diagonal scaling, also called as the Jacobi pre-conditioner is chosen, given by

\begin{equation}
\mathbf{P} = \sqrt{ \textup{diag}(\mathbf{K})^{-1} }
\end{equation}

It can be seen from \fref{fig:scaledcondk} that the scaled condition number $\mathfrak{K}$ increases with increasing degrees of freedom in case of all the approaches, viz., corrected XFEM (where blending correction with fixed area of enrichment is employed), standard XFEM ( `\textit{topological}' enrichment without blending correction) and proposed approach with different number of layers of SBFEM. But the rate of increase of the scaled condition number in case of the corrected XFEM is greater than the standard XFEM and the proposed approach. The scaled condition number of the proposed approach is lower than the standard XFEM and decreases with increasing the number of layers of the SBFEM domain. The rate of convergence of the relative error in the SIF is shown in \fref{fig:conveGriffithJint} along with conventional XFEM with only tip element enriched (\textit{'topological enrichment'}) and XFEM with two layers of elements around the tip enriched with asymptotic fields. In this example, the numerical SIF is computed using the interaction integral. It is seen that the relative error in the numerical SIF decreases with increasing mesh density and that the proposed approach yields more accurate results for same number of enrichment layers. Moreover, the relative error in the numerical SIF decreases with increasing number of layers of SBFEM domain around the crack tip for same number of elements (see \fref{fig:griffithNumlayersh}).

\begin{figure}[htpb]
\centering
\includegraphics[scale=0.6]{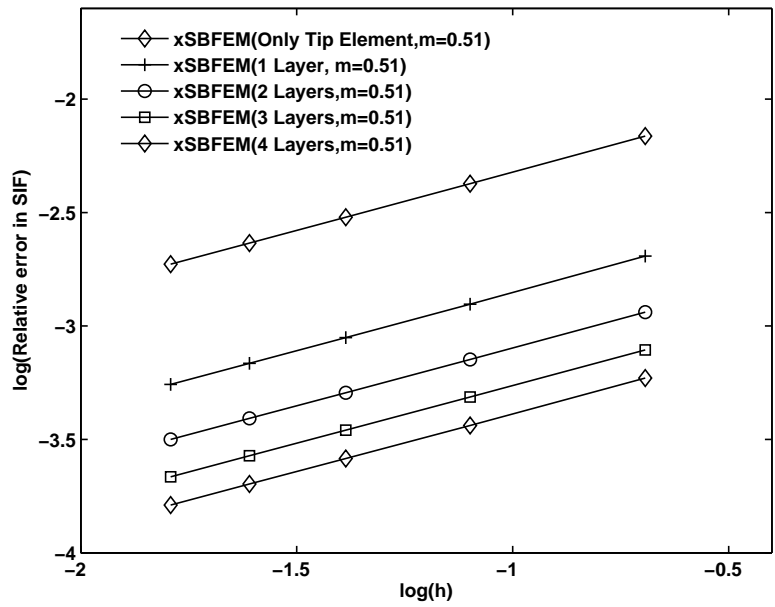}
\caption{Infinite plate with a center crack - convergence of numerical SIF with mesh density for different number of layers around the crack tip.}
\label{fig:griffithNumlayersh}
\end{figure}

\subsection{Edge crack in tension and shear loading}

\paragraph{\small Edge crack in tension}
Consider a plate with an edge crack loaded by tension $\sigma=$ 1 over the top and bottom edges. The geometry, loading and boundary conditions are show in \fref{fig:edgcrkproblem}. The reference mode I SIF is given by
\begin{equation}
K_I = F\left( \frac{a}{H} \right) \sigma \sqrt{\pi a}
\end{equation}
where $a$ is the crack length, $H$ is the plate width and $F\left( \frac{a}{H} \right)$ is an empirical function given by (for $\frac{a}{W} \le 0.6$)
\begin{equation}
F\left( \frac{a}{H} \right) = 1.12 - 0.231\left( \frac{a}{H} \right)  + 10.55 \left( \frac{a}{H} \right) ^2 -21.72 \left( \frac{a}{H} \right)^3 + 30.39 \left( \frac{a}{H} \right)^4
\end{equation}


\begin{figure}[htpb]
\centering
\scalebox{0.6}{\input{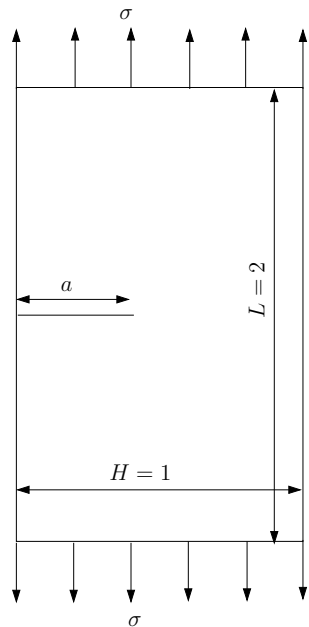}}
\caption{Plate with an edge crack under tension}
\label{fig:edgcrkproblem}
\end{figure}

\begin{figure}[htpb]
\centering
\includegraphics[scale=0.7]{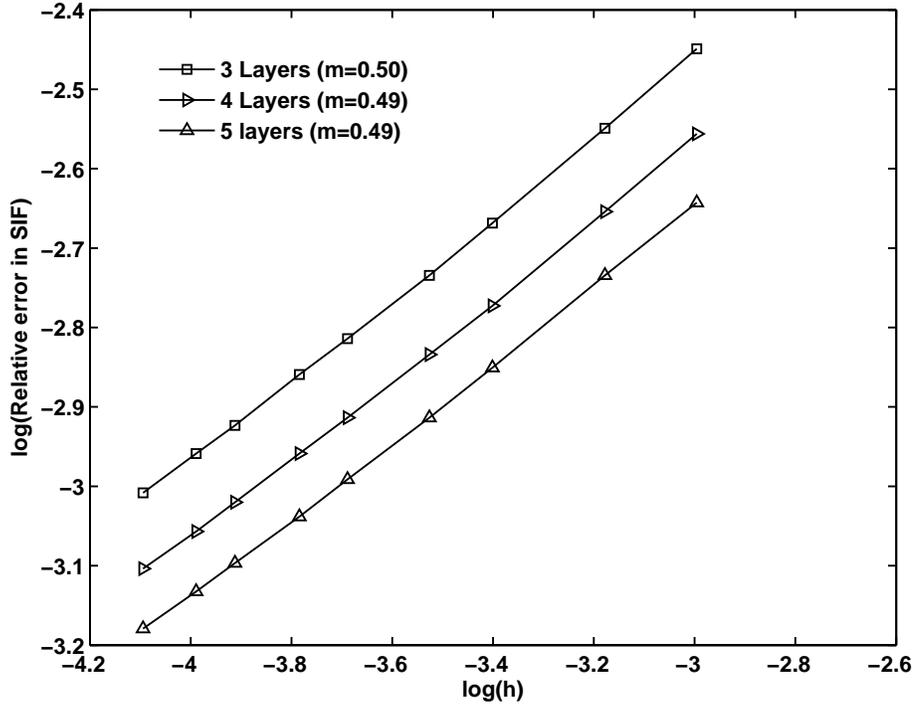}
\caption{Plate with an edge crack in tension: convergence of SIF with increasing mesh density for various number of layers around the crack tip.}
\label{fig:edgcrktenSifvsh}
\end{figure}

\begin{figure}[htpb]
\centering
\includegraphics[scale=0.7]{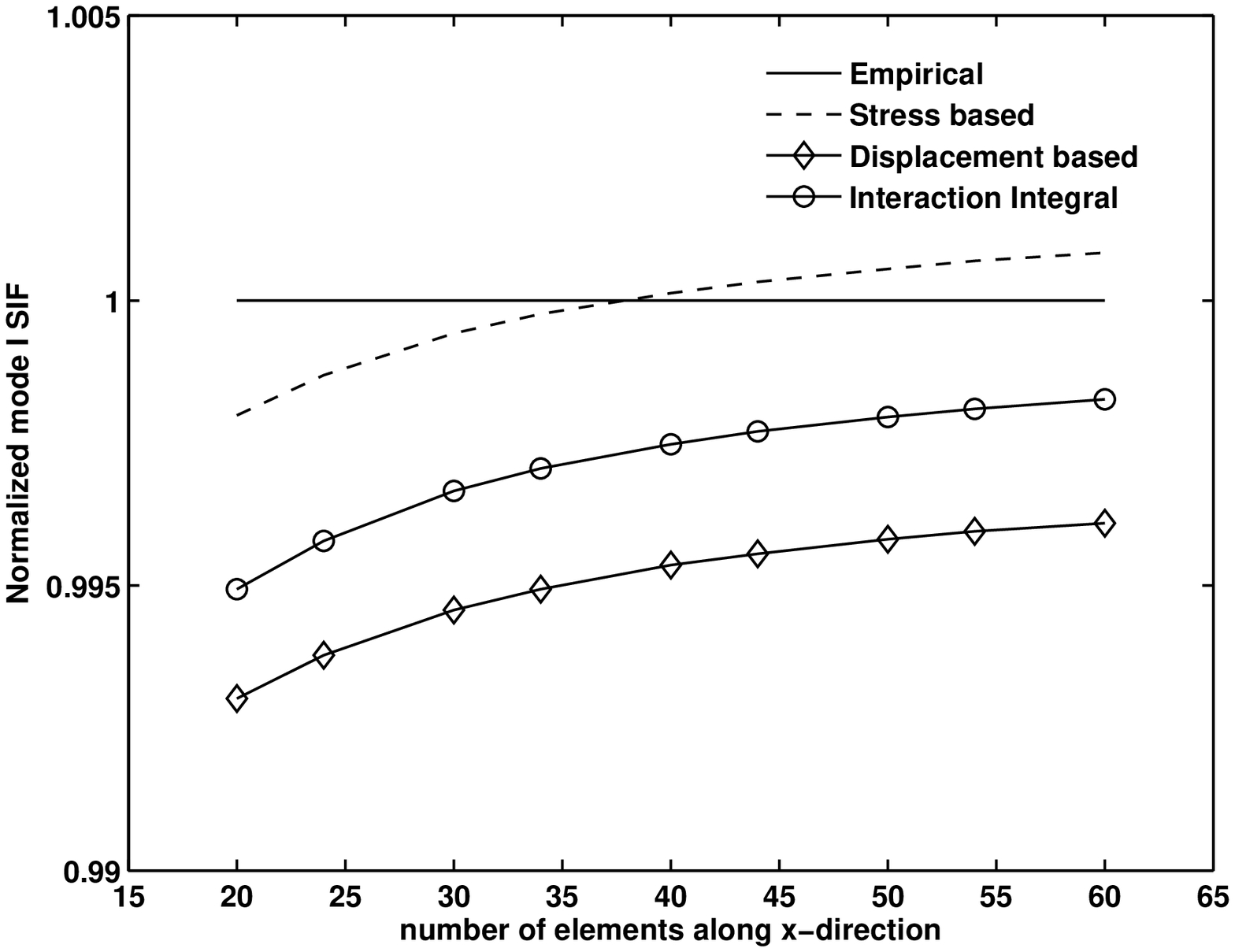}
\caption{Plate with an edge crack loaded in tension convergence results in the numerical SIF. Numerical SIF is computed using different techniques, viz., stress based, interaction integral and displacement based (see Section \ref{sifcalculation}).}
\label{fig:edgcrktension}
\end{figure}


The convergence of mode I SIF with mesh size and with the number of layers of SBFEM domain for a plate with an edge crack is shown in \fref{fig:edgcrktenSifvsh}. It can be seen that decreasing the mesh parameter $h$, the relative error in the numerical SIF decreases and that the error further decreases with increasing the size of the SBFEM domain. A similar behaviour is observed in the previous example. It can be concluded that 4 or 5 layers around the crack tip can yield accurate results. \fref{fig:edgcrktension} shows the convergence of the SIF with mesh refinement. The numerical SIF is computed by various methods, viz., stress based and displacement based method (see Section \ref{sifcalculation}) and by the interaction integral method. It can be seen that with mesh refinement, the numerical SIF approaches the empirical value for all the different approaches. In this case, 5 layers of elements around the crack tip are replaced by the SBFEM domain. 


\paragraph{\small Edge crack in shear}
In this example, consider a plate with an edge crack subjected to a shear load $\tau=$ 1 N/m$^2$ as shown in \fref{fig:edgcrkSproblem}. The material properties are assumed to be: Young's modulus, $E=$ 3$\times$10$^7$ N/m$^2$ and Poisson's ratio $\nu=$ 0.25. The reference SIFs are given as: $K_I=$ 34 N-m$^{-3/2}$ and $K_{II}=$ 4.55 N-m$^{-3/2}$. From Table \ref{tab:edgShear}, it is evident that with mesh refinement, the computed SIFs converge to the reference mixed mode SIFs. It is noted that in this example, 5 layers of elements around the crack tip are replaced with the SBFEM domain.

\begin{figure}[htpb]
\centering
\scalebox{0.6}{\input{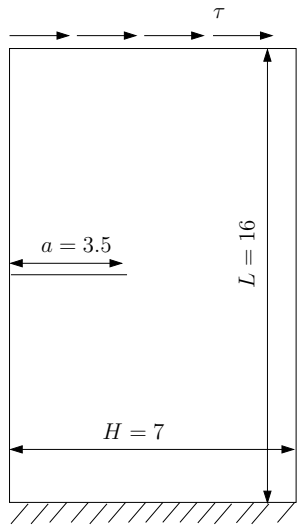}}
\caption{Plate with an edge crack under shear: geometry and boundary conditions}
\label{fig:edgcrkSproblem}
\end{figure}

\begin{table}[htbp]
\centering
\renewcommand{\arraystretch}{1.5}
\caption{Normalized mode I and mode II stress intensity factor for a plate with an edge crack with shear loading, where $K^{emp}_I=$ 34 N-m$^{-3/2}$ and $K^{emp}_{II}=$ 4.55 N-m$^{-3/2}$. }
\begin{tabular}{lrrrrrrr}
\hline
Mesh & \multicolumn{3}{c}{mode I SIF $(K_I/K^{emp}_I)$} & & \multicolumn{3}{c}{mode II SIF $(K_{II}/K^{emp}_{II})$} \\
\hline
20$\times$40	& 0.99662  & 0.99419 & 	0.99500 & & 0.99253 & 0.99301 & 0.99488 \\
30$\times$60	& 0.99871 & 0.99635 & 0.99756 & & 0.99349 & 0.99398 & 0.99587 \\
40$\times$80	& 0.99968 & 0.99737 & 0.99870 & & 0.99385 & 0.99435 & 0.99629 \\
50$\times$100 &	1.00022 & 0.99794 & 0.99933 & & 0.99402 & 0.99453 & 0.99651 \\
60$\times$120 &	1.00056 & 0.99829 & 0.99972 & & 0.99413 & 0.99464 & 0.99664 \\
\hline
\end{tabular}%
\label{tab:edgShear}%
\end{table}%



\subsection{Crack in an orthotropic body}

\paragraph{\small Center crack in an orthotropic body}
Consider a square plate $(h/w=$1$)$ with a center-crack of length $2a=$ 0.4 under uniform far field tension along the two opposite sides (see \fref{fig:orthocentcrk}. In this example, the Poisson's ratio and the shear modulus are taken as: $\nu_{12}=$ 0.03 and $G_{12}=$ 6 GPa, respectively. The Young's moduli $E_1$ and $E_2$ are computed from the following expressions~\cite{hattoridiaz2012}:
\begin{align}
E_1 &= G_{12} (\Phi + 2\nu_{12} + 1) \nonumber \\
E_2 &= \frac{E_1}{\Phi}
\end{align}
where $\Phi$ is the material parameter defined by the ratio between the Young's moduli. \fref{fig:OrthoCcrk} shows the influence of the material parameter on the normalized mode I SIF $(K_I/(\sigma\sqrt{\pi a}))$. The results are compared with those obtained in Ref.~\cite{hattoridiaz2012}. It is emphasized that in Ref.~\cite{hattoridiaz2012}, new set of enrichment functions are derived for anisotropic materials and as such no enrichment functions are required in the proposed technique. The SIF in case of the proposed method is computed using the stress definitions as outlined in Section \ref{sifcalculation}. It can be seen that the results obtained with the xSBFEM are in very good agreement with the results available in the literature. The SIF at the crack tips $A$ and $B$ are also shown in \fref{fig:OrthoCcrk} and they are similar as expected. 

\begin{figure}[htpb]
\centering
\scalebox{0.8}{\input{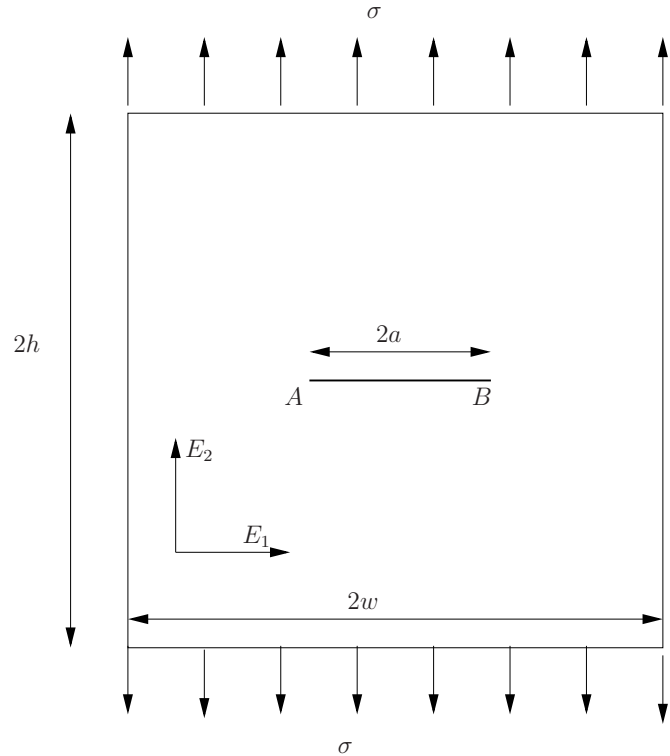}}
\caption{A square plate with a center crack under uniform far field tension. $E_1$ and $E_2$ are the Young's moduli along the $x-$ and the $y-$ directions.}
\label{fig:orthocentcrk}
\end{figure}

\begin{figure}[htpb]
\centering
\includegraphics[scale=0.65]{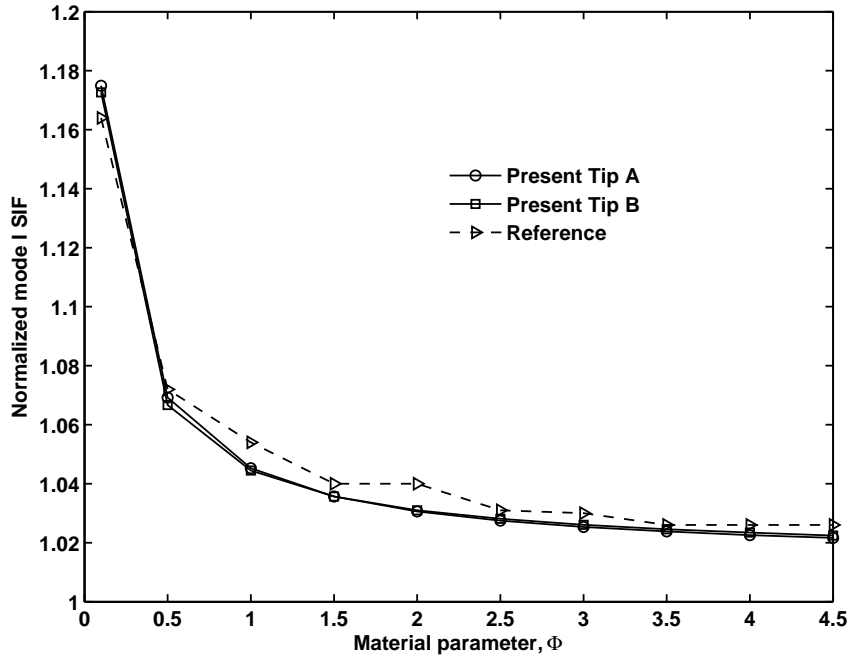}
\caption{Influence of material parameters on the normalized mode I SIF $(K_I/(\sigma \sqrt{\pi a}))$ for an orthotropic plate with a center crack.}
\label{fig:OrthoCcrk}
\end{figure}

\paragraph{\small Edge crack under uniform tension} Consider a square plate $(h/w=1)$ with a double edge-crack $(a/w=0.5)$. The plate is subjected to uniform tractions, as shown in the \fref{fig:orthoedgecrktension}. The plate is a symmetric angle ply composite laminate consisting of four graphite-epoxy laminae, with the following elastic properties: $E_1=$ 144.8 GPa, $E_2=$ 11.7 GPa, $G_{12}=$ 9.66 GPa and $\nu_{12}=$ 0.21. In this case, the influence of the fiber orientation on the SIF is studied by varying the fiber angle from $\theta=$ 0$^\circ$ to $\theta=$ 90$^\circ$. Due to symmetry, only one half of the plate is modelled. A structured quadrilateral mesh of 60$\times$120 elements is used. \fref{fig:OrthoEcrk} shows the variation of normalized mode I SIF $(K_I/(\sigma \sqrt{\pi a}))$ with the fiber orientation. It can be seen that with increasing fiber orientation, the numerical SIF increases and reaches a maximum value at $\theta=$90$^\circ$. The mode I SIF is symmetric with respect to the fiber orientation $\theta=$90$^\circ$ as expected.

\begin{figure}[htpb]
\centering
\scalebox{0.8}{\input{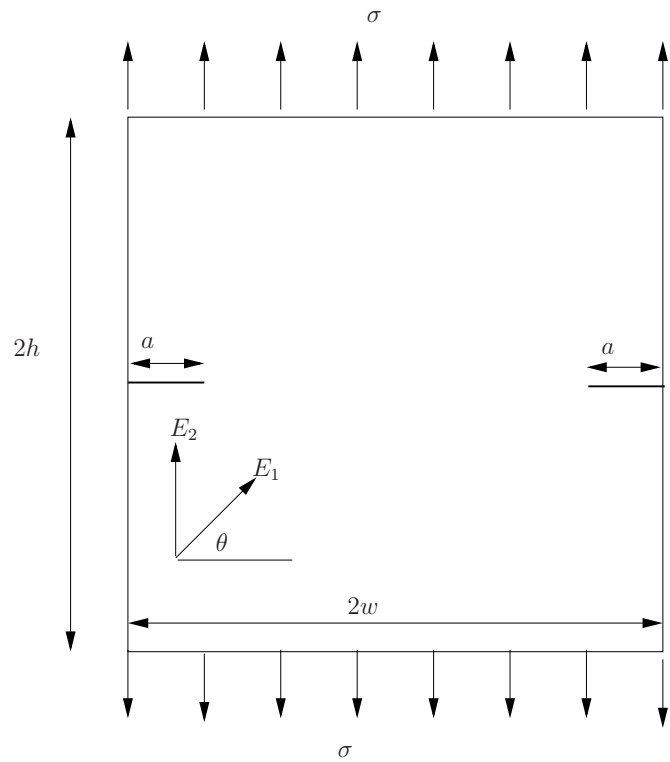}}
\caption{A square plate with a center crack under uniform far field tension. $E_1$ and $E_2$ are the Young's moduli along the $x-$ and the $y-$ directions.}
\label{fig:orthoedgecrktension}
\end{figure}

\begin{figure}[htpb]
\centering
\includegraphics[scale=0.65]{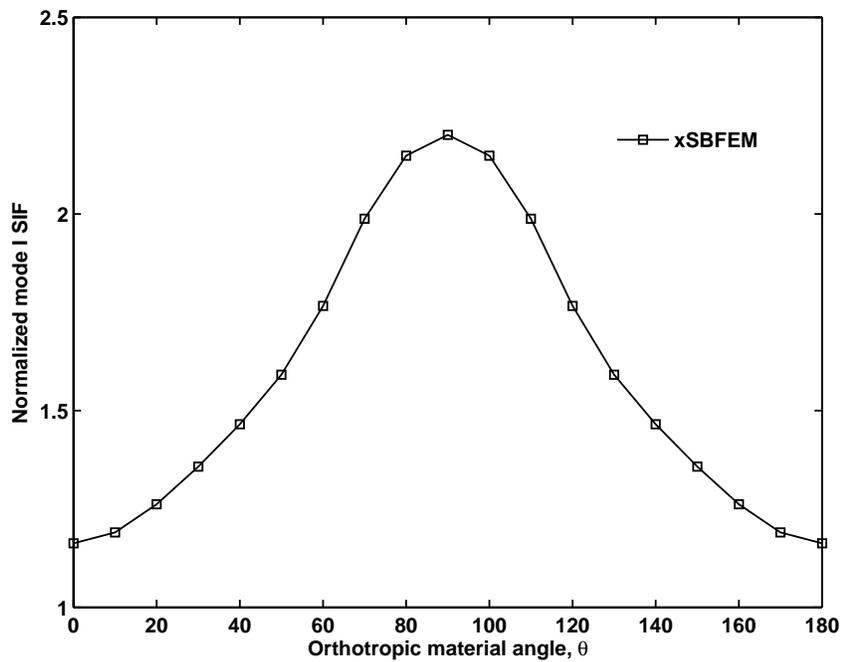}
\caption{Variation of the normalized mode I SIF $(K_I/(\sigma \sqrt{\pi a}))$ with orthotropic angle}
\label{fig:OrthoEcrk}
\end{figure}

\section{Conclusions}
In this paper, we have proposed a new method by combining the scaled boundary finite element method with the extended finite element method for problems with strong discontinuity. Instead of enriching a small region in the vicinity of the crack tip with asymptotic fields, the stiffness of the small region around the tip is computed by employing the SBFEM. With a few examples from linear elastic fracture mechanics, the effectiveness of the proposed method is illustrated. It is seen that the proposed method outperforms the conventional approach to handle discontinuities inside the domain. The key features of the proposed method without asymptotic enrichment are: (a) eliminates the need to sub-divide the element containing the crack tip; (b) suppresses the need to integrate singular functions usually appearing in the stiffness matrix, as there is no enrichment with asymptotic fields; (c) the entire process can be parallelized; (d) the scaled condition number of the proposed technique is better than the XFEM counterpart and (e) reduces the computational effort as the total number of unknowns are reduced. Although the numerical examples and the formulation presented here are for a structured mesh with bilinear elements. The proposed approach can easily be extended to other element formulations, viz., triangular, higher order quadrilateral or polygonal elements and can be applied to heterogeneous materials.

\section{Acknowledgement} Sundararajan Natarajan would like to acknowledge the financial support of the School of Civil and Environmental Engineering, The University of New South Wales for his research fellowship since September 2012.

\section{Appendix A}
\subsection{Computation of the coefficient matrices and the stiffness matrix}
For a 2-node line element with nodal coordinates ($x_{1}$, $y_{1}$), ($x_{2}$, $y_{2}$) and the shape functions along the line element are given by:
\begin{equation*}
N_{1}(\eta)= \frac{1}{2}(1-\eta); \hspace{0.25cm}
N_{2}(\eta)=  \frac{1}{2}(1+\eta).
\label{eq-twonodeshpf} 
\end{equation*}
The integrations over $\eta$ in the coefficient matrices $\ee_o, \ee_1$ and $\ee_2$ (see \Eref{eqn:coeffmat}) are performed analytically,
yielding  
\begin{eqnarray*}
\ee_0 &=& \frac{2}{3} \bq^0\left[\begin{array}{cc}
2 & 1\\
1 & 2 \end{array}\right] \nonumber \\
\ee_1 &=& \frac{1}{3} \bq^0 \left[\begin{array}{cc}
-1 & 1\\
1 & -1 \end{array}\right]+2 \bq^1\left[\begin{array}{cc}
-1 & -1\\
1 & 1
\end{array}\right] \nonumber \\
\ee_2 &=& \frac{1}{3} \bq^0 \left[\begin{array}{cc}
1 & -1\\
-1 & 1 \end{array}\right] + 4 \bq^2 \left[\begin{array}{cc}
1 & -1\\
-1 & 1 \end{array}\right]
\label{eqn:eq-2dv-2n-E012}
\end{eqnarray*}
where the following abbreviations are introduced
\begin{eqnarray*}
\bq^0 &=&  \frac{1}{4a}  \bc_1^{\rm T} \dd \bc_1 \\
\bq^1 &=& -\frac{1}{4a} \bc_2^{\rm T} \dd \bc_1 \\
\bq^2 &=& \frac{1}{4a} \bc_2^{\rm T} \dd \bc_2
\label{eqn:eq-2dv-2n-Q012} 
\end{eqnarray*} 
with 
\begin{eqnarray*}
 \bc_1 &=& \left[ \begin{array}{cc}
y_{2}-y_{1} & 0\\
0 & -(x_{2}-x_{1})\\
-(x_{2}-x_{1}) & y_{2}-y_{1}
\end{array} \right] \\
\bc_2 &=&\frac{1}{2} \left[ \begin{array}{cc}
y_{2} + y_{1} & 0 \\
0 & -(x_{2}+x_{1}) \\
-(x_{2}+x_{1}) & y_{2}+y_{1}
\end{array} \right] \\
a &=& x_{1}y_{2}-x_{2}y_{1}
\label{eqn:eq-2dv-2n-C12}
\end{eqnarray*} 

\section{Appendix B}
\begin{lstlisting}

function [ Kb, lambda, v] = getSBFEMStiffMat( coord, conn, D)

% purpose: compute the stiffness matrix for the SBFEM domain
%
% Inputs:
%	 coord - nodal coordinates
%	 conn - element connectivity along the boundary
%	 D - elasticity matrix
%
% Outputs
%	 Kb - stiffness matrix
%	 lambda, v - eigenvalues and eigenvectors
%	lambda, v  are required later for computing SIFs
%---------------------------------------------------
nd = numel(coord); id = 1:nd;
E0 = zeros(nd, nd);  E1 = zeros(nd, nd);  E2 = zeros(nd, nd);
for ie = 1:size(conn,2)
    x = coord(1,conn(1:2,ie));  y = coord(2,conn(1:2,ie));
    d = [ 2*conn(1,ie)-1 2*conn(1,ie) 2*conn(2,ie)-1 2*conn(2,ie) ];
    [ e0_ele, e1_ele, e2_ele] = EleCoeffMatrices2NodeEle(x, y, D);
    E0(d,d)=E0(d,d)+e0_ele; E1(d,d)=E1(d,d)+e1_ele; E2(d,d)=E2(d,d)+e2_ele;
end
m = E0\[E1' -eye(nd)]; Z = [m; E1*m(:,id)-E2-1d-12*E0 E1*m(:,nd+id)];
[v, d] = eig(Z);  lambda = diag(d); 
%sort eignvalues in ascending order
[~, idx] = sort(real(lambda),'ascend'); 
%rearrange eigenvalues and eigenvectors
lambda = lambda(idx(id));  v = v(:, idx(id)); 
%stiffness matrix
Kb  = real(v(nd+id, :)/v(id, :));

function [ e0, e1, e2] = EleCoeffMatrices2NodeEle( x, y, D )

% purpose: compute elemental coefficient matrices.
%
% Inputs:
%	 x,y - current element nodal coordinates
%	 D - elasticity matrix
%
% Outputs:
%	 e0,e1,e2 - SBFEM element coefficient matrices
%
%-----------------------------------------------------
a = x(1)*y(2)-x(2)*y(1);
C1 = [ y(2)-y(1) 0;   0 -(x(2)-x(1));   -(x(2)-x(1)) y(2)-y(1) ];
C2 = 0.5*[ y(2)+y(1) 0;   0 -(x(2)+x(1));   -(x(2)+x(1)) y(2)+y(1)];
Q0 = 0.25/a*(C1'*D*C1); Q1 = -0.25/a*(C2'*D*C1); Q2 = 0.25/a*(C2'*D*C2);
e0 = 2/3*[ 2*Q0 Q0; Q0 2*Q0];
e1 = 1/3*[ -Q0  Q0;  Q0 -Q0 ] + 2*[ -Q1 -Q1;   Q1 Q1 ];
e2 = 1/3*[  Q0 -Q0; -Q0  Q0 ] + 4*[  Q2 -Q2;  -Q2 Q2 ];

\end{lstlisting}

\bibliographystyle{wileyj}
\bibliography{./Biblio/myRef}
\end{document}